\documentclass[pdflatex,sn-mathphys-ay]{sn-jnl}

\usepackage{xfrac}
\usepackage{bm}
\usepackage{algorithmic}
\usepackage{natbib}
\usepackage{bibunits}
\usepackage{placeins}

\usepackage{graphicx}%
\usepackage{multirow}%
\usepackage{amsmath,amssymb,amsfonts}%
\usepackage{amsthm}%
\usepackage{mathrsfs}%
\usepackage{xcolor}%
\usepackage{textcomp}%
\usepackage{manyfoot}%
\usepackage{booktabs}%
\usepackage{algorithm}%
\usepackage{listings}%

\theoremstyle{thmstyleone}%
\newtheorem{theorem}{Theorem}
\newtheorem{proposition}[theorem]{Proposition}%
\newtheorem{lemma}[theorem]{Lemma}%

\theoremstyle{thmstyletwo}%
\newtheorem{remark}{Remark}%

\theoremstyle{thmstylethree}%

\usepackage{color}
\newcommand{\blue}{\color{black}}

\begin{document}

\title[Goemans-Williamson for Orthogonality Constraints]{Improved Approximation Algorithms for Orthogonally Constrained Problems Using Semidefinite Optimization}

\author[1]{\fnm{Ryan} \sur{Cory-Wright}}\email{r.cory-wright@imperial.ac.uk}

\author*[2]{\fnm{Jean} \sur{Pauphilet}}\email{jpauphilet@london.edu}

\affil[1]{\orgname{Imperial Business School}, \orgaddress{\street{South Kensington Campus}, \city{London}, \postcode{SW7 2AZ}, \country{United Kingdom}}}

\affil*[2]{\orgname{London Business School}, \orgaddress{\street{Regent's Park}, \city{London}, \postcode{NW1 4SA}, \country{United Kingdom}}}


\abstract{Building on the blueprint from Goemans and Williamson (1995) for the Max-Cut problem, we construct a polynomial-time approximation algorithm for orthogonally constrained quadratic optimization problems. 
First, we derive a semidefinite relaxation and propose a randomized rounding algorithm to generate feasible solutions from the relaxation. Second, we derive {constant-factor} approximation guarantees for our algorithm.
When optimizing for $m$ {orthonormal} vectors in dimension $n$, {we leverage strong duality and semidefinite complementary slackness to} show that our algorithm achieves a {$1/3$-}approximation ratio. 
{For any $m$ of the form $2^q$ for some integer $q$, we also construct an instance where the performance of our algorithm is exactly $(m+2)/(3m)$, which can be made arbitrarily close to $1/3$ by taking $m \rightarrow + \infty$, hence showing that}
our analysis is tight.
}





\maketitle

\begin{bibunit}[plainnat]

\section{Introduction}
Many important optimization problems, such as quantum {non}locality \citep{briet2011generalized}, and control theory \citep{ben2002tractable} problems feature semi-orthogonal matrices, i.e., matrices $\bm{U} \in \mathbb{R}^{n \times m}$ where $\bm{U}^\top \bm{U}=\bm{I}_m$. 
Orthogonality constraints are also related to the rank of a matrix, which models a matrix's complexity in imputation \citep{bell2007lessons}, factor analysis \citep{bertsimas2017certifiably}, and multi-task regression \citep{negahban2011estimation} settings. 

In combinatorial optimization, a major advance in the design of approximation algorithms occurred with \citet{goemans1995improved}, who proposed a $0.87856$-approximation algorithm for Max-Cut.  
Their algorithm also provides a $2/\pi$-approximation for general binary quadratic optimization (BQO) problems \citep{nesterov1998semidefinite}, and can be extended to linearly-constrained BQO problems \citep{bertsimas1998semidefinite}. 
Conceptually, \citet{goemans1995improved} established semidefinite optimization and correlated rounding at the core of approximation algorithms \citep[see][]{wolkowicz2012handbook,williamson2011design}. 

In this work, we extend the core ideas underpinning the Goemans--Williamson algorithm to quadratic semi-orthogonal 
optimization problems and provide analogous {constant-factor} guarantees on the quality of semidefinite relaxations in the semi-orthogonal setting. 

\subsection{The Original Goemans--Williamson Algorithm}\label{ssec:bqo_review}
BQO is a canonical optimization problem \citep[see][for a review of applications]{luo2010semidefinite}.
It is also an important building block for 
logically constrained optimization problems with quadratic objectives \citep[see, e.g.,][]{dong2015regularization}.
Formally, given a matrix $\bm{Q} \succeq \bm{0}$, BQO selects a vector $\bm{z}$ in $\{-1, 1\}^{m}$ that solves
\begin{align}\label{prob:maxcut2}
    \max_{\bm{z} \in \{-1, 1\}^{m}} \quad & {\sum_{i,j} Q_{i,j} z_i z_j =\max_{\bm{z} \in \{-1, 1\}^{m}}} \langle \bm{Q}, \bm{z}\bm{z}^\top\rangle,
\end{align}
where $\langle \cdot, \cdot \rangle$ denotes the Frobenius inner product between matrices.
Problem \eqref{prob:maxcut2} is NP-hard and often challenging to solve to optimality when $m \geq 100$ \citep{rehfeldt2023faster}. Accordingly, a popular approach for obtaining near-optimal solutions is to sample from a distribution parameterized by the solution of \eqref{prob:maxcut2}'s convex relaxation. Specifically, we can reformulate \eqref{prob:maxcut2} in terms of the rank-one matrix $\bm{Z} = \bm{z}\bm{z}^\top$. Then, 
a valid relaxation of Problem \eqref{prob:maxcut2} is given by
\begin{align}\label{prob:maxcut3}
    \max_{\bm{Z} \in \mathcal{S}_+^{m}} \quad & \langle \bm{Q},  \bm{Z}\rangle \ \text{s.t.} \ \mathrm{diag}(\bm{Z})=\bm{e},
\end{align}
which would be an exact reformulation with the additional (non-convex) constraint $\operatorname{rank}(\bm{Z}) = 1$.
Probabilistically speaking, \eqref{prob:maxcut3} is a device for constructing a pseudodistribution over $\bm{z} \in \{-1, 1\}^{m}$ \citep{d2003relaxations}. This suggests  
sampling vectors from a distribution with second moment $\bm{Z}^\star$ 
and rounding to restore feasibility, as proposed by \citet{goemans1995improved} for Max-Cut and described in Algorithm~\ref{alg:gwmethod}. 
The overall idea of Algorithm \ref{alg:gwmethod} is that 
the projection step (i.e., taking the coordinate-wise sign of $\bm{y}$) partially preserves the second moment of the distribution of $\bm{y}$, $\mathbb{E}[{\bm y} {\bm y}^\top ] = \bm{Z}^\star$.
Precisely, we have $\mathbb{E}\left[\hat{\bm{z}}\hat{\bm{z}}^\top \right] \succeq \frac{2}{\pi} \bm{Z}^\star$ \citep[see][]{nesterov1998semidefinite, bertsimas1998semidefinite}.
\begin{algorithm}
\caption{The Goemans--Williamson rounding algorithm for Problem \eqref{prob:maxcut2}}\label{alg:gwmethod}
\begin{algorithmic}\normalsize
\STATE Compute $\bm{Z}^\star$ {an optimal} solution of \eqref{prob:maxcut3}
\STATE Sample $\bm{y} \sim \mathcal{N}(\bm{0}, \bm{Z}^\star)$
\STATE Construct $\hat{\bm{z}} \in \{-1, 1\}^{m}: \hat{z}_i := \mathrm{sign}(y_i)$
\RETURN $\hat{\bm{z}}$ a solution to Problem \eqref{prob:maxcut2}
\end{algorithmic}
\end{algorithm}

\subsection{Orthogonally Constrained Quadratic Optimization}
In this paper, we consider a family of orthogonally constrained quadratic problems that subsumes binary quadratic optimization. Formally, we search for $m$ orthogonal vectors $\bm{u}_i \in \mathbb{R}^n$ which solve
\begin{align}
    \max_{\bm{u}_i \in \mathbb{R}^{n},\ i \in [m]} \quad \sum_{i,j = 1}^m \bm{u}_i^\top \bm{A}^{(i,j)} \bm{u}_j \mbox{ s.t. } \bm{u}_i^\top \bm{u}_j = \delta_{i,j}, \quad \forall i,j \in [m],
\end{align}
where $\bm{A}$ is an {$nm \times nm$ positive} semidefinite matrix {(in short, $\bm{A} \in \mathcal{S}^{nm}_+$)} with block matrices $\bm{A}^{(i,j)} \in \mathbb{R}^{n \times n}$, {and} $\delta_{i,j}=1$ if $i=j$ and $0$ otherwise. 
We require $n \geq m$. By introducing the semi-orthogonal matrix $\bm{U} \in \mathbb{R}^{n \times m}$ whose columns are the vectors $\bm{u}_i \in \mathbb{R}^n$, we can write our problem as
\begin{align}\label{prob:orthogonallyconstrained}
    \max_{\bm{U} \in \mathbb{R}^{n \times m}} \quad \operatorname{vec}(\bm{U})^\top \bm{A} \operatorname{vec}(\bm{U}) \quad & \text{s.t.} \quad  \bm{U}^\top \bm{U}=\bm{I}_m,
\end{align}
where the $\operatorname{vec}(\cdot)$ operator stacks the column{s} of $\bm{U}$ together into a single vector. 

The similarities between Problems \eqref{prob:orthogonallyconstrained} and \eqref{prob:maxcut2} are striking: For example, we can {reduce} any BQO instance \eqref{prob:maxcut2} {to} a special case of Problem \eqref{prob:orthogonallyconstrained}.
{Our reduction (see Appendix \ref{sec:a.bqo.equiv}) not only shows that Problem \eqref{prob:orthogonallyconstrained} is NP-hard \citep[as also proved in][Theorem 3.1]{lai2025stiefel} but is also approximation-preserving. 
Therefore, inheriting the inapproximability results of BQO, Problem~\eqref{prob:orthogonallyconstrained} cannot be approximated in polynomial time within a factor of ${2/\pi}+\varepsilon$ unless P=NP
\citep[][Theorem I.3]{briet2015tight}. 
}

In addition to its applications in clustering, quantum nonlocality, or generalized trust-region problems \citep[see][and references therein]{burer2023strengthened}, Problem \eqref{prob:orthogonallyconstrained} appears as a relevant substructure for mixed-projection formulations of rank-constrained quadratic optimization problems \citep{bertsimas2020mixed}.

In this paper, inspired by the Goemans--Williamson algorithm for BQO, we develop a relax-then-round strategy with a {$1/3$}-approximation ratio for Problem \eqref{prob:orthogonallyconstrained}. 

\subsection{Related Work} \label{ssec:orthogonallyconstrained}
Our work is most closely related to \citet{burer2023strengthened}, who develop a hierarchy of semidefinite relaxations for Problem~\eqref{prob:orthogonallyconstrained}. To numerically evaluate the tightness of their relaxations, they apply several `feasible rounding procedures' but do not provide any theoretical approximation ratios. In contrast, we develop a randomized rounding procedure and show that it achieves a multiplicative factor guarantee, which is independent of {both} the ambient dimension $n$ and {the number of vectors $m$}. As a non-convex quadratic optimization problem, 
Problem~\eqref{prob:orthogonallyconstrained} can also be solved to optimality via global solvers such as \verb|Gurobi| or \verb|BARON|. However, the scalability of these global solvers is currently limited. 

\subparagraph{Special Cases} A larger body of work considers a special case of Problem \eqref{prob:orthogonallyconstrained}, where the matrix $\bm{A}$ is block-diagonal, namely $\bm{A}^{(i,j)}= \bm{0}$ if $i \neq j$. In this case, Problem \eqref{prob:orthogonallyconstrained} reduces to
\begin{align}\label{prob:orthogonallyconstrained_blockcase}
    \max_{\bm{U} \in \mathbb{R}^{n \times m}} \sum_{i \in [m]}\quad \bm{u}_i^\top \bm{A}^{(i,i)} \bm{u}_i \quad & \text{s.t.} \quad  \bm{U}^\top \bm{U}=\bm{I}_m,
\end{align}
which is referred to as the sum of heterogeneous quadratic forms or the heterogeneous PCA problem. Indeed, when all the matrices $\bm{A}^{(i,i)} $ are equal, we recover the Principal Component Analysis (PCA) problem. \citet[Section 5]{bolla1998extrema} solve Problem \eqref{prob:orthogonallyconstrained_blockcase} in polynomial time via linear algebra techniques when the matrices $\bm{A}^{(i,i)}$ are diagonal or commute with each other. For general matrices, \cite{gilman2022semidefinite} further tailor the semidefinite relaxations of \citet{burer2023strengthened}. Although tighter, their relaxations are not always tight. For some instances, they even obtain optimality gaps exceeding $100\%$.
We are not aware of any approximation algorithms with guarantees {specifically} for general (non-diagonal) instances of Problem \eqref{prob:orthogonallyconstrained_blockcase}.

\subparagraph{Approximation Algorithms} To our knowledge, {many of the} existing approximation algorithms {apply} to optimization problems with different orthogonality structures. \citet{briet2010positive} propose an approximation algorithm for problems of the form
\begin{align}
    \max_{\bm{U} \in \mathbb{R}^{n \times m}} \sum_{i,j \in [m]} A_{i,j} \bm{u}_i^\top \bm{u}_j 
    \ \text{s.t.} \ \bm{u}_i^\top \bm{u}_i=1 \ \forall i \in [m],
\end{align}
which also subsumes BQO (for $n=1$), but does not enforce orthogonality between the columns of $\bm{U}$. They devise a relax-and-round strategy analogous to Goemans--Williamson that achieves an approximation ratio of $2/\pi+\Theta(1/n)$. 

A second line of work \citep{nemirovski2007sums,man2009improved} proposes $\Omega(1/\log(n+m))$-approximation algorithms for quadratic optimization problems over matrices $\bm{U}$ that satisfy $\bm{U}^\top \bm{U} \preceq \bm{I}_m$, {i.e., whose largest singular value is at most one}. This constraint does not ensure that the columns of $\bm{U}$ are orthogonal. 
Nonetheless, \citet{nemirovski2007sums} shows that, in several {special} cases such as the orthogonal Procrustes or quadratic assignment problems (but not in the case of Problem~\eqref{prob:orthogonallyconstrained}), orthogonality constraints can be relaxed into $\bm{U}^\top \bm{U} \preceq \bm{I}_m$  without loss of optimality. 
{Our algorithm differs in that it generates matrices $\bm{U}$ with orthogonal columns, while they only guarantee $\sigma_{\max}(\bm{U}) \leq 1$.
In addition, their approximation ratio vanishes both with $n$ and $m$, while ours is dimension-independent.} 

Finally, \citet{bandeira2016approximating} study approximation algorithms for problems of the form
\begin{align} \label{eqn:bandeira}
    \max_{\bm{U}_i \in \mathbb{R}^{n \times m}, i=1,\dots,k} \: \sum_{i, j \in [k]}\langle \bm{A}^{(i,j)}, \bm{U}_i^\top \bm{U}_{j} \rangle \ \text{s.t.} \ \bm{U}_i^\top \bm{U}_i=\bm{I}_{m} \ \forall i \in [k].
\end{align}
Unfortunately, Problem \eqref{eqn:bandeira} is not equivalent to \eqref{prob:orthogonallyconstrained}. 
In particular, heterogeneous PCA is a special case of our Problem \eqref{prob:orthogonallyconstrained} but cannot be cast in the form \eqref{eqn:bandeira}.
There are two key differences in the objective function of \eqref{eqn:bandeira}: it involves the \emph{inner} products between columns of \emph{different} semi-orthogonal matrices $\bm{U}_i$, $\bm{U}_{i'}$ for $i\neq i'$. 
On the other hand, the objective in \eqref{prob:orthogonallyconstrained} depends on \emph{outer} products between columns of the \emph{same} matrix $\bm{U}$. In particular, \citet{bandeira2016approximating} can restore feasibility for each $\bm{U}_i$, $i=1,\dots,k$ in \eqref{eqn:bandeira} separately, while the columns of $\bm{U}$ in \eqref{prob:orthogonallyconstrained} need to be orthogonalized together. {Finally, \citet{bandeira2016approximating}'s proof technique involves taking a singular value decomposition of a random matrix and arguing that the singular vectors are distributed according to the Haar measure and independent of the singular values. Unfortunately, our setting \eqref{prob:orthogonallyconstrained} is more general than \eqref{eqn:bandeira}, and the singular vectors in the random matrices studied here are not distributed according to the Haar measure.}

{
However, both our Problem \eqref{prob:orthogonallyconstrained} with $\bm{A} \succeq \bm{0}$ and Problem \eqref{eqn:bandeira} are special cases of the non-commutative Grothendieck problem. To clarify this connection, let us reformulate \eqref{prob:orthogonallyconstrained} as a bilinear optimization problem, using a trick analogous to that of \citet[][Section 5.3]{naor2014efficient}.
\begin{lemma} Assume that $\bm{A} \succeq \bm{0}$. Then, Problem \eqref{prob:orthogonallyconstrained} achieves the same objective value as
    \begin{align}\label{prob:orthogonallyconstrained.bilinear}
    \max_{\bm{U}, \bm{V} \in \mathbb{R}^{n \times m}} \quad \langle \bm{A}, \operatorname{vec}(\bm{U}) \operatorname{vec}(\bm{V})^\top\rangle\quad & \text{s.t.} \quad  \bm{U}^\top \bm{U}=\bm{V}^\top \bm{V}=\bm{I}_m.
\end{align}
\end{lemma}
\begin{proof} Any feasible solution $\bm{U}$ to Problem \eqref{prob:orthogonallyconstrained} defines a feasible solution $(\bm{U}, \bm{U})$ to Problem \eqref{prob:orthogonallyconstrained.bilinear} with the same objective value, so \eqref{prob:orthogonallyconstrained} $\leq$ \eqref{prob:orthogonallyconstrained.bilinear}. Conversely, given a feasible solution to \eqref{prob:orthogonallyconstrained.bilinear}, we can apply Cauchy-Schwarz because $\bm{A} \succeq \bm{0}$ and get
\begin{align*}
  \langle \bm{A}, \operatorname{vec}(\bm{U}) \operatorname{vec}(\bm{V})^\top\rangle   \leq \| \bm{A}^{1/2} \operatorname{vec}(\bm{U}) \|_2 \| \bm{A}^{1/2} \operatorname{vec}(\bm{V}) \|_2 \leq \eqref{prob:orthogonallyconstrained},
\end{align*}
which concludes the proof.
\end{proof}
By embedding $\bm{U}$ and $\bm{V}$ as the first $m$ columns of $n \times n$ orthogonal matrices, $\widetilde{\bm{U}}$ and $\widetilde{\bm{V}}$, and appropriately padding the matrix $\bm{A}$ with zeros, Problem \eqref{prob:orthogonallyconstrained.bilinear} can be expressed as an instance of the non-commutative Grothendieck problem studied, for example, in \citet{naor2014efficient} and \citet{briet2015tight}. \citet{naor2014efficient} propose a $1/(2\sqrt{2})$-factor approximation algorithm for this problem \citep[][Figure 3]{naor2014efficient}. Although our approximation ratio is marginally weaker, our algorithm relies on a more compact semidefinite relaxation (involving one $nm \times nm$ semidefinite matrix vs. $2n^2 \times 2 n^2$), and is therefore more tractable, especially when $n \gg m$. It also applies directly to real-valued variables (rather than their complex-valued rounding step) and recovers the original Goemans--Williamson algorithm for BQO instances. 
}

\subsection{Contributions and Structure}\label{ssec:contributions}
Our main contribution is the development of a Goemans--Williamson sampling algorithm for the class of semi-orthogonal problems \eqref{prob:orthogonallyconstrained}.

In Section \ref{sec:gw_matrix}, we derive a semidefinite relaxation and propose a sampling procedure. {Our algorithm achieves a $1/3$-approximation guarantee (Theorem \ref{thm:13}) for Problem \eqref{prob:orthogonallyconstrained}. Our analysis is tight in the sense that for any $\epsilon > 0$, we can construct an instance such that the approximation ratio of our algorithm is at most $1/3 + \epsilon$ (Proposition \ref{prop:ce.13}). While our approximation ratio does not depend on $n,m$, we can also derive tighter dimension-dependent approximation ratios (Proposition \ref{prop:2pim} and Theorem \ref{thm:chisquared.light}).
The proofs of our main results are presented in Section \ref{sec:proof}. In particular, our analysis relies heavily on properties of the dual of the semidefinite relaxation.}
To better judge the quality of our approximations, we develop two simpler algorithms in Section~\ref{sec:benchmark}: uniform sampling and a stronger benchmark inspired by PCA. We show that they achieve $1/nm$ and $ 1/m^2$ approximation ratios, respectively, which are dominated by our semidefinite {relax-and-project} procedure. Then, we evaluate the performance of our approach numerically in Section~\ref{sec:numerics}.

\subsection{Notations} \label{ssec:notation}
We let lowercase boldfaced characters such as $\bm{x}$ denote vectors and uppercase boldfaced characters such as $\bm{X}$ denote matrices. We denote by $\bm{e}$ the vector of all ones. 
We let $[n]$ denote the set of running indices $\{1, ..., n\}$. The cone of $n \times n$ symmetric (resp. positive {semi}definite) matrices is denoted by $\mathcal{S}^n$ (resp. $\mathcal{S}^n_+$). 
For a matrix $\bm{X} \in \mathbb{R}^{n \times m}$, we let 
$\bm{x}_i$ denote its $i$th column. 
We let $\operatorname{vec}(\bm{X}): \mathbb{R}^{n \times m} \rightarrow \mathbb{R}^{nm}$ denote the vectorization operator which maps matrices to vectors by stacking columns. 
For a matrix $\bm{W}$, we may describe it as a block matrix composed of equally sized blocks and denote the $(i,i')$ block by $\bm{W}^{(i,i')}$. The dimension of each block will be clear from the context. 
{For any matrices $\bm{A}, \bm{B}$, we denote their Kronecker product $\bm{A} \otimes \bm{B}$. We will extensively use the identity $\operatorname{vec}(\bm{ABC}) = (\bm{C}^\top \otimes \bm{A}) \operatorname{vec}(\bm{B})$ \citep[see, e.g.,][]{henderson1981vec}.}
Finally, we let $\mathcal{N}(\bm{0}, \bm{\Sigma})$ denote a centered multivariate normal distribution with covariance matrix $\bm{\Sigma}$.

\section{A Goemans--Williamson Approach}\label{sec:gw_matrix}
In this section, we propose a new Goemans--Williamson-type approach for semi-orthogonal quadratic optimization problems. {First}, we review a semidefinite relaxation for semi-orthogonal quadratic optimization in Section~\ref{ssec:gw_matrix.relax}. {Second,} in Section~\ref{ssec:gw_matrix.sampling}, we propose a randomized rounding scheme to generate feasible solutions. {Finally,} we present a multiplicative approximation ratio of $1/3$ and a dimension-dependent approximation ratio for our algorithm in Section \ref{ssec:gw_matrix.perf}.

\subsection{A Shor Relaxation}\label{ssec:gw_matrix.relax}
As reviewed in Section \ref{ssec:orthogonallyconstrained}, \citet[][Section 2.2]{burer2023strengthened} derive the following semidefinite relaxation for Problem~\eqref{prob:orthogonallyconstrained}:
\begin{equation}\label{prob:orth_relax_shor}
\begin{aligned}
\max_{\bm{W} \in \mathcal{S}^{mn}_+} \:  \langle \bm{A}, \bm{W} \rangle \quad \text{s.t.} \: & \sum_{i \in [m]}\bm{W}^{(i,i)} \preceq \bm{I}_n, \ \operatorname{tr}\left( \bm{W}^{(j,j')} \right)= \delta_{j,j'}, \forall j, j' \in [m],
\end{aligned}
\end{equation}
where the matrix $\bm{W}$ encodes the outer-product of $\operatorname{vec}(\bm{U})$ with itself, and the trace constraints on the blocks of $\bm{W}$ stem from the columns of $\bm{U}$ having unit norm and being pairwise orthogonal. 

Similarly to the semidefinite relaxation of \eqref{prob:maxcut2}, imposing the constraint that $\bm{W}$ is rank-one in \eqref{prob:orth_relax_shor} would result in an exact reformulation of \eqref{prob:orthogonallyconstrained}. Accordingly, the relaxation \eqref{prob:orth_relax_shor} is tight whenever some optimal solution is rank-one. 
However, the optimal solutions to \eqref{prob:orth_relax_shor} are often high-rank {(the case $m=1$ is one of the special cases where this semidefinite relaxation is tight).}
Actually, it follows from manipulating the Barvinok--Pataki bound \citep{barvinok2001remark,pataki1998rank} that there exists some optimal solution to Problem \eqref{prob:orth_relax_shor} with rank at most $n+m$.\footnote{After introducing a slack matrix $\bm{S}$ to write the semidefinite inequality constraint $\bm{S}=\bm{I}_n-\sum_{i \in [m]}\bm{W}^{(i,i)}, \bm{S} \succeq \bm{0}$, we have $m(m+1)/2+n(n+1)/2$ scalar inequalities. Thus, \citet[Theorem 2.2]{pataki1998rank} states that there exists some optimal solution $(\bm{W}, \bm{S})$ with $\mathrm{rank}(\bm{W})(\mathrm{rank}(\bm{W})+1)/2+\mathrm{rank}(\bm{S})(\mathrm{rank}(\bm{S})+1)/2\leq m(m+1)/2+n(n+1)/2$. Since $\mathrm{rank}(\bm{S}) \geq 0$, it implies $\mathrm{rank}(\bm{W})\leq n+m$.} However, not all optimal solutions obey this bound; thus, we do not use this observation in our analysis.
An interesting question is how to generate a high-quality feasible solution to \eqref{prob:orthogonallyconstrained}, with a provable approximation ratio, by leveraging a solution of 
\eqref{prob:orth_relax_shor}, which is the focus of the rest of the section. 

Note that Problem \eqref{prob:orth_relax_shor}, which is essentially a Shor relaxation \citep{shor1987quadratic}, corresponds to the `DiagSum' relaxation of \citet{burer2023strengthened}. They also derive an even stronger relaxation, which they call a `Kronecker' relaxation. 
We do not explicitly analyze their Kronecker relaxation here because it is significantly less tractable \citep[as reported in][Table 1]{burer2023strengthened}. 

\subsection{A Relax-and-Project Procedure} \label{ssec:gw_matrix.sampling}
We propose a randomized rounding scheme to generate high-quality feasible solutions to \eqref{prob:orthogonallyconstrained} from an optimal solution to \eqref{prob:orth_relax_shor}.

Our algorithm involves three main steps: First, we solve \eqref{prob:orth_relax_shor} and obtain a semidefinite matrix $\bm{W}^\star$. Second, using $\bm{W}^\star$, we sample an $n \times m$ matrix $\bm{G}$ such that $\operatorname{vec}(\bm{G})$ follows a normal distribution with mean $\bm{0}_{nm}$ and covariance matrix $\bm{W}^\star$. Third, from the matrix $\bm{G}$, we generate a feasible solution to \eqref{prob:orthogonallyconstrained} {by projecting $\bm{G}$ onto the set of semi-orthogonal matrices}. 
Specifically, we compute a singular value decomposition of $\bm{G}$, $\bm{G} = \bm{U\Sigma V}^\top$ and define $\bm{Q} := \bm{UV}^\top$.
We summarize our procedure in Algorithm~\ref{alg:matrixgwalgorithm2}.

In Algorithm \ref{alg:matrixgwalgorithm2}, we can sample $\operatorname{vec}(\bm{G}) \sim \mathcal{N}(\bm{0}_{nm}, \bm{W}^\star)$ even when $\bm{W}^\star$ is rank-deficient via the following construction---which will also be relevant for our theoretical analysis. Denoting $r = \operatorname{rank}(\bm{W}^\star)$, we first construct a Cholesky decomposition of $\bm{W}^\star$: $\bm{W}^\star = \sum_{k \in [r]} \operatorname{vec}(\bm{B}_k)\operatorname{vec}(\bm{B}_k)^\top$ with $\bm{B}_k \in \mathbb{R}^{n \times m}$. Then, we sample $\operatorname{vec}(\bm{G}) = \sum_{k \in [r]} \operatorname{vec}(\bm{B}_k) z_k$ with $\bm{z} \sim \mathcal{N}(\bm{0}_r, \bm{I}_r)$. 
This procedure ensures that $\operatorname{vec}(\bm{G}) \in \operatorname{span}(\bm{W}^\star)$ almost surely. 
In particular, {if our semidefinite relaxation is tight (e.g., when $m=1$) and $\bm{W}^\star$ has rank 1, then $\bm{Q}$ is optimal almost surely (our rounding is exact)}.
\begin{algorithm}
\caption{A relax-and-project algorithm for Problem \eqref{prob:orthogonallyconstrained}}
\label{alg:matrixgwalgorithm2}
\begin{algorithmic}
    \STATE Compute $\bm{W}^\star$ {an optimal} solution of \eqref{prob:orth_relax_shor}
    \STATE Sample $\bm{G}$ according to $\operatorname{vec}(\bm{G}) \sim \mathcal{N}(\bm{0}_{nm}, \bm{W}^\star)$
    \STATE Compute the SVD of $\bm{G}$, $\bm{G} = \bm{U\Sigma V}^\top$
    \STATE Construct $\bm{Q} =\bm{U}\bm{V}^\top$
    \RETURN {the} semi-orthogonal matrix $\bm{Q}$ 
\end{algorithmic}
\end{algorithm}

{By properties of the Gaussian distribution, we can explicitly control some second and fourth moments of the random matrix $\bm G$, which will be crucial in deriving our constant-factor approximation guarantees. 
 \begin{lemma} \label{lemma:moments} Consider $\bm{W} \in \mathcal{S}_+^{nm}$ a feasible solution to \eqref{prob:orth_relax_shor}. The random matrix $\bm{G}$ generated by Algorithm \ref{alg:matrixgwalgorithm2} satisfies
\begin{align*}
\mathbb{E}[\bm{G}^\top \bm{G}] = \bm{I}_m, 
\quad \mbox{ and } \quad \mathbb{E}[\bm{GG}^\top] \preceq \bm{I}_n.
\end{align*}
Moreover it satisfies
\begin{align*}
\mathbb{E}[(\bm{G}^\top \bm{G})^2] \preceq 3\bm{I}_m, 
\quad \mbox{ and } \quad  \mathbb{E}[(\bm{GG}^\top)^2] \preceq 3\bm{I}_n.
\end{align*}
\end{lemma}
The equations on second moments in Lemma \ref{lemma:moments} stem from the constraints imposed in our semidefinite relaxation. Indeed, denoting 
$\bm{g}_j,j=1,\dots,m$ the columns of $\bm{G}$, we have that the $(j,j')$ block of $\bm{W}$, $\bm{W}^{(j,j')}$, corresponds to $\mathbb{E}[\bm{g}_j \bm{g}_{j'}^\top]$. In addition, we can express the matrices $\bm{G}^\top \bm{G}$ and $\bm{G}\bm{G}^\top$ as$[\bm{g}_j^\top \bm{g}_{j'}]_{j,j'}$ and $\sum_{j\in [m]} \bm{g}_j \bm{g}_j^\top$ respectively.  
The bounds on the fourth moments are a generalization of the known equality $\mathbb{E}[z^4] =3$ for $z\sim\mathcal{N}(0,1)$ to matrix Gaussian series. We defer a formal proof of Lemma \ref{lemma:moments} to Section \ref{ssec:a.moments}.}

Similar to the original algorithm of \citet{goemans1995improved}, the intuition behind Algorithm \ref{alg:matrixgwalgorithm2} is that the sampled matrix $\bm{G}$ achieves an average performance equal to the relaxation value ($\mathbb{E}[\langle \bm{A}, \operatorname{vec}(\bm{G})\operatorname{vec}(\bm{G})^\top\rangle]=\langle \bm{A}, \bm{W}^\star\rangle$) and is feasible on average ($\mathbb{E}[\bm{G}^\top \bm{G}]=\bm{I}_m$). Therefore, the {average} objective value of the feasible solution $\bm{Q}$ should not be too different from that of $\bm{G}$, as we now theoretically study.

\subsection{A Constant-Factor Guarantee} \label{ssec:gw_matrix.perf}
{In this section, we study the theoretical} performance of Algorithm \ref{alg:matrixgwalgorithm2} in the case where the objective matrix $\bm{A}$ in \eqref{prob:orthogonallyconstrained} is positive semidefinite.

{
Our main result is a constant-factor guarantee for Algorithm~\ref{alg:matrixgwalgorithm2}: 
\begin{theorem}\label{thm:13} Assume that $\bm{A} \succeq \bm{0}$. Let $\bm{W}^\star$ denote an optimal solution to the semidefinite relaxation \eqref{prob:orth_relax_shor}. The random matrix $\bm{Q} \in \mathbb{R}^{n \times m}$ generated by Algorithm~\ref{alg:matrixgwalgorithm2} satisfies the inequality
\begin{align*}
\mathbb{E}\left[\operatorname{vec}(\bm{Q})^\top \bm{A} \operatorname{vec}(\bm{Q}) \right] & \geq  \dfrac{1}{3} \langle \bm{A}, \bm{W}^\star \rangle.
\end{align*}
\end{theorem}
We defer the proof of Theorem \ref{thm:13} to Section \ref{sec:proof}.
Theorem~\ref{thm:13} provides a worst-case approximation ratio of $1/3$ for Algorithm~\ref{alg:matrixgwalgorithm2}, which is independent of the ambient dimension $n$ as well as the number of vectors $m$. Although the actual performance of our algorithm can be much stronger (as we demonstrate numerically in Section \ref{sec:numerics}), we can show that our worst-case analysis is tight. 

\begin{proposition}\label{prop:ce.13}
    Consider an integer $m$ of the form $m = 2^q$ for some $q \geq 1$. There exists an integer $n \geq m$, a positive semidefinite matrix $\bm{A}$, and an optimal solution $\bm{W}$ to the semidefinite relaxation \eqref{prob:orth_relax_shor}, such that the solutions generated by Algorithm \ref{alg:matrixgwalgorithm2} satisfy
    \begin{align*}
        \operatorname{vec}(\bm{Q})^\top \bm{A} \operatorname{vec}(\bm{Q}) = \dfrac{m+2}{3m} \langle \bm{A},\bm{W} \rangle \quad a.s.
     \end{align*}
\end{proposition}
For any $\epsilon > 0$, applying Proposition~\ref{prop:ce.13} (proof in Section \ref{ssec:a.proof.ce13}) for a sufficiently large $m$ shows the existence of an instance on which the approximation ratio of Algorithm \ref{alg:matrixgwalgorithm2} is at most $1/3 +\epsilon$, hence showing that, from a worst-case perspective, the $1/3$-approximation ratio of Theorem \ref{thm:13} is essentially tight. 

However, it does not rule out the existence of algorithms with tighter guarantees or the possibility to derive tighter constants for fixed values of $n$ and $m$. In particular, we can show the following (proof deferred to Section \ref{ssec:a.proof.2pim}): 
\begin{proposition} \label{prop:2pim}
Assume that $\bm{A} \succeq \bm{0}$. Let $\bm{W}^\star$ denote an optimal solution to the semidefinite relaxation \eqref{prob:orth_relax_shor}. The random matrix $\bm{Q} \in \mathbb{R}^{n \times m}$ generated by Algorithm~\ref{alg:matrixgwalgorithm2} satisfies the inequality
\begin{align*}
\mathbb{E}\left[\operatorname{vec}(\bm{Q})^\top \bm{A} \operatorname{vec}(\bm{Q}) \right] & \geq  \dfrac{2}{\pi m} \langle \bm{A}, \bm{W}^\star \rangle.
\end{align*}
\end{proposition}
For $m = 1$, Proposition \ref{prop:2pim} provides a $2/\pi \approx 0.636$ approximation factor, which is strictly better than $1/3$. However, for $m \geq 2$, the guarantee of Proposition~\ref{prop:2pim} is weaker than that of Theorem \ref{thm:13}.  

\begin{remark} In Theorem~\ref{thm:chisquared.light}, we derive an even tighter approximation guarantee, which decreases monotonically with $n$ and $m$. In particular, for $m=1$, our improved constant is at least equal to 0.6804 (in the limit $n \rightarrow \infty$). For $m=2$, it is at least equal to $0.3402$. However, it fails to improve upon the 1/3 guarantee of Theorem \ref{thm:13} for $m > 2$.
\end{remark}
}

{
\section{Proof of the Main Results} 
\label{sec:proof}
The proof of Theorem \ref{thm:13} and Proposition \ref{prop:2pim} rely on strong duality and complementarity slackness results for the semidefinite relaxation \eqref{prob:orth_relax_shor}. We present these intermediate technical results in Section \ref{ssec:proof.dualitylemma}. We then present the proof of Theorem \ref{thm:13} and Proposition \ref{prop:2pim} in Sections \ref{ssec:proof.13} and \ref{ssec:a.proof.2pim} respectively. 

\subsection{Semidefinite relaxation and duality} \label{ssec:proof.dualitylemma}
We first derive the dual of Problem \eqref{prob:orth_relax_shor}.
\begin{lemma} \label{lemma:relax.dual}
Assume that $\bm{A}\in\mathcal{S}_+^{mn}$. The dual of \eqref{prob:orth_relax_shor} is equal to
\begin{equation}\label{eq:dual}
\begin{aligned}
\min_{\bm{X}\in\mathcal{S}^n_+,\ \bm{\Lambda}\in\mathcal{S}^m}\quad & \operatorname{tr}(\bm{X})+\operatorname{tr}(\bm{\Lambda})\\
\text{s.t.}\quad & \bm{I}_m\otimes \bm{X} + \bm{\Lambda}\otimes \bm{I}_n \succeq \bm{A}.
\end{aligned}
\end{equation}
and strong duality holds. 
Furthermore, we can take $\bm{\Lambda} \succeq \bm{0}$ in \eqref{eq:dual} without loss of optimality.
\end{lemma}
The fact that we can take $\bm{\Lambda} \succeq \bm{0}$ will be crucial in our proof. Indeed, with this additional assumption, for any positive semidefinite matrices $\bm{M}_1,\bm{M}_2$, the relationship $\bm{M}_1 \succeq \bm{M}_2$ implies $\operatorname{tr}(\bm{\Lambda}\bm{M}_1) \geq \operatorname{tr}(\bm{\Lambda}\bm{M}_2)$, which would not be true if $\bm{\Lambda} \not\succeq \bm{0}$.
\begin{proof}
We introduce a semidefinite dual variable $\bm{X}\in\mathcal{S}^n_+$ for the primal semidefinite constraint
$\sum_{i \in [m]} \bm{W}^{(i,i)}\preceq \bm{I}_n$, and scalars $\lambda_{j,j'}$ for the trace equalities
$\operatorname{tr}(\bm{W}^{(j,j')})=\delta_{j,j'}$, for all $(j,j') \in \{1,\dots,m\}^2$. Collectively, the variables $(\lambda_{j,j'})$ form a symmetric matrix
$\bm{\Lambda}\in\mathcal{S}^m$. The Lagrangian is defined 
\begin{align*}
\mathcal{L}(\bm{W},\bm{X},\bm{\Lambda})
&=\langle \bm{A},\bm{W}\rangle
+\Big\langle \bm{X},\bm{I}_n-\sum_{i \in [m]} \bm{W}^{(i,i)}\Big\rangle
+\sum_{j,j' \in [m]} \lambda_{j,j'}\big(\delta_{j,j'}-\operatorname{tr}(\bm{W}^{(j,j')})\big).
\end{align*}

Denoting $\{ \bm{E}_{j,j'} \}_{j,j'=1,\dots,m}$ the canonical basis of $\mathbb{R}^{m \times m}$, we can concisely rewrite
\begin{align*}
\sum_{i\in [m]} \langle \bm{X},\bm{W}^{(i,i)}\rangle = \sum_{i \in [m]} \langle \bm{E}_{i,i} \otimes \bm{X},\bm{W} \rangle 
= \langle \bm{I}_{m} \otimes \bm{X},\bm{W} \rangle,
\end{align*}
and
\begin{align*}
\sum_{j,j' \in [m]} \lambda_{j,j'} \operatorname{tr}(\bm{W}^{(j,j')})= \sum_{j,j' \in [m]} \lambda_{j,j'} \langle \bm{E}_{j,j'} \otimes \bm{I}_n, \bm{W} \rangle 
= \langle \bm{\Lambda} \otimes \bm{I}_n, \bm{W} \rangle
\end{align*}
Taking the supremum of the Lagrangian over $\bm{W} \in \mathcal{S}^{nm}_+$ then yields the dual of the form \eqref{eq:dual}. 

Observe that 
$(\bm{X},\bm{\Lambda}) = (\bm{I}_n,\lambda_{\max}(\bm{A}) \bm{I}_m)$ is a strictly feasible solution to Problem \eqref{eq:dual}. Therefore, Slater's conditions are satisfied and strong duality holds \citep[][Theorem 3.1]{vandenberghe1996semidefinite}.

For the second part of the proof, consider a feasible dual solution $\bm{X} \in\mathcal{S}^n_+, \bm{\Lambda} \in\mathcal{S}^m$. The semidefinite constraint in \eqref{eq:dual} implies
\begin{align*}
\bm{I}_m\otimes \bm{X} + \bm{\Lambda}\otimes \bm{I}_n \succeq \bm{A} \succeq \bm{0}.
\end{align*}
Take a unit vector $\bm{u} \in \mathbb{R}^m$ and $\bm{x} \in \mathbb{R}^n$. We have
\begin{align*}
&(\bm{u} \otimes \bm{x})^\top (\bm{I}_m\otimes \bm{X} + \bm{\Lambda}\otimes \bm{I}_n) (\bm{u} \otimes \bm{x}) = \bm{x}^\top \bm{X} \bm{x} + (\bm{u}^\top \bm{\Lambda} \bm{u}) \| \bm{x} \|^2 \geq 0 \end{align*}
which implies $\bm{X} \succeq - \lambda_{\min}(\bm{\Lambda}) \bm{I}_n$. Set $\gamma := \max(0,- \lambda_{\min}(\bm{\Lambda})) \geq 0$. We have $\bm{X} \succeq \gamma \bm{I}_n$. 

Define $\bm{X}^\star = \bm{X} - \gamma \bm{I}_n$, $\bm{\Lambda}^\star = \bm{\Lambda} + \gamma \bm{I}_m$. 
Then, $\bm{X}^\star \succeq \bm{0}$ and $\bm{I}_m\otimes \bm{X}^\star + \bm{\Lambda}^\star\otimes \bm{I}_n  = \bm{I}_m\otimes \bm{X} + \bm{\Lambda}\otimes \bm{I}_n$. So, $(\bm{X}^\star, \bm{\Lambda}^\star)$ is feasible for \eqref{eq:dual}. By construction, $\bm{\Lambda}^\star\succeq \bm{0}$. In terms of objective value, $\operatorname{tr}(\bm{X}^\star)+\operatorname{tr}(\bm{\Lambda}^\star) = \operatorname{tr}(\bm{X})+\operatorname{tr}(\bm{\Lambda}) - \gamma (n-m) \leq \operatorname{tr}(\bm{X})+\operatorname{tr}(\bm{\Lambda})$ because $n \geq m$. So, $(\bm{X}^\star, \bm{\Lambda}^\star)$ achieves a lower objective value.
\end{proof}

The second technical lemma uses complementarity slackness of a primal-dual pair to simplify bilinear expressions of the form $\operatorname{vec}(\bm{Q})^\top \bm{A} \operatorname{vec}(\bm{G})$.
\begin{lemma}\label{lemma:bilin} Consider $\bm{W}^\star$ an optimal solution to the semidefinite relaxation \eqref{prob:orth_relax_shor} and $(\bm{X}^\star, \bm{\Lambda}^\star) \in \mathcal{S}^n_+ \times \mathcal{S}^m_+$ an optimal solution to the dual problem \eqref{eq:dual}.

Consider a random matrix $\bm{G} \in \mathbb{R}^{n \times m}$ generated as in Algorithm \ref{alg:matrixgwalgorithm2}: $\operatorname{vec}(\bm{G}) \sim \mathcal{N}(\bm{0}, \bm{W}^\star)$ and denote $\bm{Q}$ its projection onto the set of semi-orthogonal matrices. We have 
\begin{align}
    \operatorname{tr}(\bm{\Lambda}^\star \mathbb{E}[\bm{G}^\top\bm{G}]) = \operatorname{tr}(\bm{\Lambda}^\star) \quad \mbox{ and } \quad \operatorname{tr}(\bm{X}^\star \mathbb{E}[\bm{GG}^\top]) = \operatorname{tr}(\bm{X}^\star). \label{eqn:zero.tr}
\end{align}
Furthermore, the following relationships hold almost surely:
\begin{align}
 \operatorname{vec}(\bm{Q})^\top \bm{A} \operatorname{vec}(\bm{G})
 &= \operatorname{tr}(\bm{X}^\star (\bm{G} \bm{G}^\top)^{1/2}) + \operatorname{tr}(\bm{\Lambda}^\star (\bm{G}^\top \bm{G})^{1/2}) \label{eqn:bilin.qg} \\ 
 &\geq \dfrac{1}{\| \bm{G} \|_F} \operatorname{vec}(\bm{G})^\top \bm{A} \operatorname{vec}(\bm{G}) \label{eqn:bilin.frob}.
\end{align}
\end{lemma}
\begin{proof}
The first equality in \eqref{eqn:zero.tr} follows directly from the fact that $\mathbb{E}[\bm{G}^\top\bm{G}] = \bm{I}_m$ (Lemma \ref{lemma:moments}). 

For the second equality, $\mathbb{E}[\bm{G}\bm{G}^\top] \preceq \bm{I}_n$ (Lemma \ref{lemma:moments}) only provides an inequality. To obtain an equality, we leverage the fact that $(\bm{W}^\star, \bm{X}^\star, \bm{\Lambda}^\star)$ are primal-dual optimal. In particular, they satisfy  $\langle \bm{X}^\star,\bm{I}_n-\sum_{i\in [m]} \bm{W}^{\star (i,i)}\rangle = 0$ by complementarity slackness. Given that $\mathbb{E}[\bm{G}\bm{G}^\top] = \sum_{i\in [m]} \bm{W}^{\star (i,i)}$, we obtain $$\operatorname{tr}(\bm{X}^\star \mathbb{E}[\bm{GG}^\top]) = \operatorname{tr}(\bm{X}^\star).$$

Equations \eqref{eqn:bilin.qg}-\eqref{eqn:bilin.frob} also follow from complementarity slackness. From $\langle \bm{W}^\star, \bm{I}_m\otimes \bm{X}^\star + \bm{\Lambda}^\star\otimes \bm{I}_n - \bm{A} \rangle = 0$, we get that 
\begin{align*}
    (\bm{I}_m\otimes \bm{X}^\star + \bm{\Lambda}^\star\otimes \bm{I}_n - \bm{A}) \bm{W}^\star = \bm{0},
\end{align*}
because both matrices are positive semidefinite.
In other words, the range of $\bm{W}^\star$ is included in the kernel of the matrix $(\bm{I}_m\otimes \bm{X}^\star + \bm{\Lambda}^\star\otimes \bm{I}_n - \bm{A})$. Since, by construction, $\operatorname{vec}(\bm{G}) \in \operatorname{span}(\bm{W}^\star)$ almost surely (as discussed in Section \ref{ssec:gw_matrix.sampling}), we must have
\begin{align*}
    (\bm{I}_m\otimes \bm{X}^\star + \bm{\Lambda}^\star\otimes \bm{I}_n - \bm{A}) \operatorname{vec}(\bm{G}) &= \bm{0},
\end{align*}
i.e., 
\begin{align}
     \bm{A} \operatorname{vec}(\bm{G})  =  (\bm{I}_m\otimes \bm{X}^\star + \bm{\Lambda}^\star\otimes \bm{I}_n) \operatorname{vec}(\bm{G}). \label{eqn:dual.quad.g}
\end{align}
Hence, 
\begin{align*}
 \operatorname{vec}(\bm{Q})^\top \bm{A} \operatorname{vec}(\bm{G})
 &= \operatorname{vec}(\bm{Q})^\top (\bm{I}_m\otimes \bm{X}^\star + \bm{\Lambda}^\star\otimes \bm{I}_n) \operatorname{vec}(\bm{G}) 
 = \operatorname{tr}(\bm{X}^\star \bm{Q} \bm{G}^\top) + \operatorname{tr}(\bm{\Lambda}^\star \bm{Q}^\top \bm{G}).
\end{align*}
Equation \eqref{eqn:bilin.qg} follows immediately after observing that $\bm{Q} = (\bm{G} \bm{G}^\top)^{-1/2} \bm{G} = \bm{G} (\bm{G}^\top \bm{G})^{-1/2}$, so 
$\bm{Q}\bm{G}^\top = (\bm{G} \bm{G}^\top)^{1/2}$ and $\bm{Q}^\top \bm{G} = (\bm{G}^\top \bm{G})^{1/2}$.

Finally, to obtain \eqref{eqn:bilin.frob}, we apply the bound $\bm{\Sigma} \succeq \bm{\Sigma}^2 / \| \bm{\Sigma} \|_F$ and get $\bm{Q}\bm{G}^\top = (\bm{G} \bm{G}^\top)^{1/2} \succeq (\bm{G} \bm{G}^\top) /  \| \bm{G} \|_F$ and $\bm{Q}^\top \bm{G} = (\bm{G}^\top \bm{G})^{1/2} \succeq (\bm{G}^\top \bm{G}) /  \| \bm{G} \|_F$. Because $\bm{X}^\star \succeq \bm{0}$ and $\bm{\Lambda}^\star \succeq \bm{0}$, these semidefinite bounds lead to
\begin{align*}
 \operatorname{vec}(\bm{Q})^\top \bm{A} \operatorname{vec}(\bm{G})
 &\geq \dfrac{1}{\| \bm{G} \|_F} \operatorname{tr}(\bm{X}^\star \bm{G} \bm{G}^\top) + \dfrac{1}{\| \bm{G} \|_F}\operatorname{tr}(\bm{\Lambda}^\star \bm{G}^\top \bm{G}) = \dfrac{1}{\| \bm{G} \|_F} \operatorname{vec}(\bm{G})^\top \bm{A} \operatorname{vec}(\bm{G}),
\end{align*}
where the equality follows from applying \eqref{eqn:dual.quad.g} again. 
\end{proof}

\subsection{Proof of Theorem \ref{thm:13}} \label{ssec:proof.13}
The proof of Theorem \ref{thm:13} leverages the bounds on the second and fourth moments of $\bm{G}^\top \bm{G}$ and $\bm{G} \bm{G}^\top$ in Lemma \ref{lemma:moments}. However, the average performance of Algorithm \ref{alg:matrixgwalgorithm2} involves the second moment of $\operatorname{vec}(\bm{Q})$, $\mathbb{E}[\operatorname{vec}(\bm{Q}) \operatorname{vec}(\bm{Q})^\top]$. Instead of trying to control this matrix directly, the proof follows three steps: (i) We apply Cauchy-Schwarz to lower bound $\mathbb{E}[ \operatorname{vec}(\bm{Q})^\top \bm{A} \operatorname{vec}(\bm{Q})]$ by $\mathbb{E}[ \operatorname{vec}(\bm{Q})^\top \bm{A} \operatorname{vec}(\bm{G})]$; (ii) We use duality results from Section \ref{ssec:proof.dualitylemma} (in particular, Equation \eqref{eqn:bilin.qg} from Lemma \ref{lemma:bilin}) to show that $\mathbb{E}[ \operatorname{vec}(\bm{Q})^\top \bm{A} \operatorname{vec}(\bm{G})]$ only depends on $\mathbb{E}[\bm{Q}^\top \bm{G}]$  and $\mathbb{E}[\bm{Q}\bm{G}^\top ]$ ; (iii) We leverage concentration inequalities for matrix Gaussian series from Lemma \ref{lemma:moments} combined with the following technical lemma (proof deferred to Section \ref{ssec:a.12}).
\begin{lemma}\label{lemma:matrixgaussian.12} Consider a random positive semidefinite matrix $\bm{S} \in \mathcal{S}^d_+$ and a fixed semidefinite matrix $\bm{M} \in \mathcal{S}^d_+$ such that $\operatorname{tr}(\bm{M}\mathbb{E}[\bm{S}]) = \operatorname{tr}(\bm{M})$ and $\operatorname{tr}(\bm{M}\mathbb{E}[\bm{S}^{2}])\leq c \operatorname{tr}(\bm{M})$ for some constant $c > 0$. Then, $\operatorname{tr}(\bm{M}\mathbb{E}[\bm{S}^{1/2}]) \geq (1/\sqrt{c})\operatorname{tr}(\bm{M})$.
\end{lemma}

\begin{proof}[Proof of Theorem \ref{thm:13}]
Given that $\bm{A} \succeq \bm{0}$, by Cauchy-Schwarz, we have 
\begin{align*}
    \mathbb{E}[\operatorname{vec}(\bm{Q})^\top \bm{A} \operatorname{vec}(\bm{Q})] \geq \dfrac{\mathbb{E}[\operatorname{vec}(\bm{Q})^\top \bm{A} \operatorname{vec}(\bm{G})]^2}{\mathbb{E}[\operatorname{vec}(\bm{G})^\top \bm{A} \operatorname{vec}(\bm{G})]}.
\end{align*}
Let us introduce optimal dual solutions to the semidefinite relaxation, $(\bm{X}^\star, \bm{\Lambda}^\star) \in \mathcal{S}_+^n \times \mathcal{S}_+^m$, whose existence is guaranteed by Lemma \ref{lemma:relax.dual}. By strong duality, we have $\langle \bm{A}, \bm{W}^\star \rangle = \operatorname{tr}(\bm{X}^\star) + \operatorname{tr}(\bm{\Lambda}^\star)$. Furthermore, Equation \eqref{eqn:bilin.qg} in Lemma \ref{lemma:bilin} yields
\begin{align*}
\mathbb{E}[\operatorname{vec}(\bm{Q})^\top \bm{A} \operatorname{vec}(\bm{G})] = \operatorname{tr}(\bm{X}^\star \mathbb{E}[(\bm{G} \bm{G}^\top)^{1/2}]) + \operatorname{tr}(\bm{\Lambda}^\star \mathbb{E}[(\bm{G}^\top \bm{G})^{1/2}]).
\end{align*}

The random matrix $\bm{G} \bm{G}^\top$ satisfies $\operatorname{tr}(\bm{X}^\star \mathbb{E}[\bm{G} \bm{G}^\top]) = \operatorname{tr}(\bm{X}^\star)$ (Equation \eqref{eqn:zero.tr} in Lemma \ref{lemma:bilin}) and $\operatorname{tr}(\bm{X}^\star \mathbb{E}[(\bm{G} \bm{G}^\top)^2] \leq 3 \operatorname{tr}(\bm{X}^\star)$ (see Lemma \ref{lemma:moments}), so, Lemma \ref{lemma:matrixgaussian.12} implies that $$\operatorname{tr}(\bm{X}^\star \mathbb{E}[(\bm{G} \bm{G}^\top)^{1/2}]) \geq \operatorname{tr}(\bm{X}^\star) / \sqrt{3}.$$ Similarly, with $\bm{G}^\top \bm{G}$, we get 
$\operatorname{tr}(\bm{\Lambda}^\star \mathbb{E}[(\bm{G}^\top \bm{G})^{1/2}]) \geq \operatorname{tr}(\bm{\Lambda}^\star) / \sqrt{3}$.  

As a result, we have
\begin{align*}
&\mathbb{E}[\operatorname{vec}(\bm{Q})^\top \bm{A} \operatorname{vec}(\bm{G})] \geq \dfrac{1}{\sqrt{3}} \left(\operatorname{tr}(\bm{X}^\star) + \operatorname{tr}(\bm{\Lambda}^\star) \right) = \dfrac{1}{\sqrt{3}} \langle \bm{A},\bm{W}^\star \rangle, \\
\mbox{ and }
& \mathbb{E}[\operatorname{vec}(\bm{Q})^\top \bm{A} \operatorname{vec}(\bm{Q})] \geq \dfrac{1}{3}  \langle \bm{A},\bm{W}^\star \rangle.
\end{align*}
\end{proof}

\subsection{Proof of Proposition \ref{prop:2pim}} \label{ssec:a.proof.2pim}
The proof of Proposition \ref{prop:2pim} adopts a similar strategy as that of Theorem \ref{thm:13} but uses Equation \eqref{eqn:bilin.frob} from Lemma \ref{lemma:bilin} instead of Equation \eqref{eqn:bilin.qg} to lower bound $\operatorname{vec}(\bm{Q})^\top \bm{A} \operatorname{vec}(\bm{Q})$.
\begin{proof}[Proof of Proposition \ref{prop:2pim}] Given that $\bm{A} \succeq \bm{0}$, by Cauchy-Schwarz, we have 
\begin{align*}
    \operatorname{vec}(\bm{Q})^\top \bm{A} \operatorname{vec}(\bm{Q}) \geq \dfrac{(\operatorname{vec}(\bm{Q})^\top \bm{A} \operatorname{vec}(\bm{G}))^2}{\operatorname{vec}(\bm{G})^\top \bm{A} \operatorname{vec}(\bm{G})}.
\end{align*}
From Equation \eqref{eqn:bilin.frob} in Lemma \ref{lemma:bilin}, we have 
\begin{align*}
    \operatorname{vec}(\bm{Q})^\top \bm{A} \operatorname{vec}(\bm{G}) \geq \dfrac{1}{\|\bm{G}\|_F}{\operatorname{vec}(\bm{G})^\top \bm{A} \operatorname{vec}(\bm{G})},
\end{align*}
so 
\begin{align*}
    \operatorname{vec}(\bm{Q})^\top \bm{A} \operatorname{vec}(\bm{Q}) \geq 
    \dfrac{1}{\|\bm{G}\|_F^2}{\operatorname{vec}(\bm{G})^\top \bm{A} \operatorname{vec}(\bm{G})}.
\end{align*}
Taking expectation yields
\begin{align*}
    \mathbb{E}[\operatorname{vec}(\bm{Q})^\top \bm{A} \operatorname{vec}(\bm{Q})] \geq \mathbb{E}\left[ \dfrac{\operatorname{vec}(\bm{G})^\top \bm{A} \operatorname{vec}(\bm{G})}{\|\bm{G}\|_F^2} \right] = \langle \bm{A}, \mathbb{E}\left[ \dfrac{\operatorname{vec}(\bm{G}) \operatorname{vec}(\bm{G})^\top}{\|\bm{G}\|_F^2} \right] \rangle
\end{align*}
To conclude, we show that  
\begin{align*}
    \mathbb{E}\left[ \dfrac{\operatorname{vec}(\bm{G}) \operatorname{vec}(\bm{G})^\top}{\|\bm{G}\|_F^2} \right] \succeq \dfrac{2}{\pi m} \bm{W}^\star.
\end{align*}

Consider a fixed vector $\bm{v} \in \mathbb{R}^{nm}$. From Cauchy-Schwarz, we have that for, any random variables $A \geq 0, B > 0$,
$\mathbb{E}[A/B] \geq \mathbb{E}[\sqrt{A}]^2 / \mathbb{E}[B]$. Applied to $A = (\bm{v}^\top \operatorname{vec}(\bm{G}))^2$ and $B = \|\bm{G}\|_F^2$, we get 
\begin{align*}
\mathbb{E}\left[ \dfrac{ (\bm{v}^\top\operatorname{vec}(\bm{G}))^2}{\|\bm{G}\|_F^2} \right] \geq \dfrac{\mathbb{E}\left[ |\bm{v}^\top\operatorname{vec}(\bm{G})| \right]^2}{\mathbb{E}\left[\|\bm{G}\|_F^2\right]} 
\end{align*}
For the numerator, $\bm{v}^\top\operatorname{vec}(\bm{G}) \sim \mathcal{N}(0, \bm{v}^\top \bm{W}^\star \bm{v})$ so $\mathbb{E}\left[ |\bm{v}^\top\operatorname{vec}(\bm{G})| \right] = \sqrt{\tfrac{2}{\pi} \bm{v}^\top \bm{W}^\star \bm{v}}$. For the denominator, $\mathbb{E}\left[\| \bm{G} \|_F^2 \right] = \mathbb{E}\left[ \operatorname{tr}(\bm{G}^\top \bm{G}) \right] = m$, which concludes the proof. 
\end{proof}

\begin{remark}{The proof of Proposition \ref{prop:2pim} shows that 
\begin{align*}
    \mathbb{E}[\operatorname{vec}(\bm{Q})^\top \bm{A} \operatorname{vec}(\bm{Q})] \geq \left\langle \bm{A}, \mathbb{E}\left[ \dfrac{\operatorname{vec}(\bm{G}) \operatorname{vec}(\bm{G})^\top}{\|\bm{G}\|_F^2} \right] \right\rangle.
\end{align*} However, we would like to emphasize that, while this scalar inequality holds, 
we do not necessarily have the semidefinite relationship 
\begin{align*}
    \mathbb{E}[\operatorname{vec}(\bm{Q})\operatorname{vec}(\bm{Q})^\top] \succeq \mathbb{E}\left[ \dfrac{\operatorname{vec}(\bm{G}) \operatorname{vec}(\bm{G})^\top}{\|\bm{G}\|_F^2} \right].
\end{align*}
Indeed,}
consider the counterexample with $n=4, m=2$ and 
\begin{align*}
    \bm{W}^{\star} = \begin{pmatrix} \bm{A} & \bm{B} \\ \bm{B}^\top & \bm{C} \end{pmatrix}, 
\end{align*}
with 
\begin{align*} \blue
    \bm{A} & =\mathrm{Diag}\begin{pmatrix} 0.019 \\ 0.132 \\ 0.268 \\ 0.581\end{pmatrix}, \quad
     \bm{B} = \mathrm{Diag}\begin{pmatrix} 0.016 \\ 0.193\\ 0.192 \\ -0.401\end{pmatrix}, \quad  \bm{C}=\mathrm{Diag}\begin{pmatrix} 0.035 \\ 0.294 \\ 0.154 \\ 0.517 \end{pmatrix}.
\end{align*}

For each seed $s$, we sample $N=25,000$ matrices $(\bm{G}, \bm{Q})$ as in Algorithm \ref{alg:matrixgwalgorithm2} and compute $\lambda_s = \lambda_{\min}\left(\widehat{\mathbb{E}}\left[\operatorname{vec}(\bm{Q})\operatorname{vec}(\bm{Q})^\top-\frac{\operatorname{vec}(\bm{G})\operatorname{vec}(\bm{G})^\top}{ \|\bm{G}\|_F^2} \right]\right)$. 
We generate $B=100$ bootstrapped samples from the same $N$ pairs of matrices $(\bm{G}, \bm{Q})$ and compute the corresponding minimum eigenvalue $\lambda_{s,b}$. The average variance $\sigma_s^2 = \dfrac{1}{B} \sum_{b=1}^B (\lambda_{s,b} - \lambda_s)^2$ captures the within-seed error in estimating $\lambda_{\min}$ due to using a finite number of observations $N$ to approximate the expectation.
We repeat this process over $s=50$ seeds to evaluate the estimation error due to using these particular $N$ samples. An unbiased estimate of the smallest eigenvalue is given by $\hat{\lambda} = \dfrac{1}{S}\sum_{s=1}^S \lambda_s$, with standard error $\dfrac{1}{S^2}\sum_{s=1}^S (\lambda_s-\hat{\lambda})^2 + \dfrac{1}{S^2}\sum_{s=1}^S \sigma_s^2$, and we construct a $95\%$ confidence interval. We find $$\lambda_{\min}\left(\mathbb{E}\left[\operatorname{vec}(\bm{Q})\operatorname{vec}(\bm{Q})^\top-\frac{\operatorname{vec}(\bm{G})\operatorname{vec}(\bm{G})^\top}{\|\bm{G}\|_F^2} \right]\right)= -0.0136 \pm 0.0003 <0.$$ 
\end{remark}
}

\section{Benchmark: Uniform Sampling and Deflation}\label{sec:benchmark}
To appreciate the strength of our approximation ratios for Algorithm \ref{alg:matrixgwalgorithm2}, we analyze the performance of two baselines: {uniform sampling (Section \ref{ssec:benchmark.uniform}) and a deflation heuristic inspired by PCA (Section \ref{ssec:benchmark.deflation}).
The results of this section are summarized in Table \ref{tab:guarantee.summary}.}
\begin{table}[h!]
    \centering
    \caption{Summary of our guarantees for Algorithm  \ref{alg:matrixgwalgorithm2} and two benchmarks.}\vspace{2mm}
    \label{tab:guarantee.summary}
    \begin{tabular}{cccc}
        Name & Algorithm \ref{alg:matrixgwalgorithm2} & Uniform & Deflation \\ \midrule
        Ratio\ & $\dfrac{1}{3}$ & $\dfrac{1}{nm}$ & $\dfrac{1}{m^2}$\\[1em]
        Source & Theorem \ref{thm:13} & Proposition \ref{prop:haar_1} & Proposition \ref{prop:m2eigenapprox}
    \end{tabular}
\end{table}

\subsection{Uniform Sampling} \label{ssec:benchmark.uniform}
We consider a naive baseline where we draw $\bm{Q}$ uniformly from the set of semi-orthogonal matrices. Note that this is analogous to generating i.i.d. Bernoulli vectors in BQO, which achieves a $1/2$ approximation ratio in the Max-Cut case {\citep[e.g.,][Theorem 5.3]{williamson2011design}}. 

{In practice, to sample $\bm{U}$ uniformly from the set of semi-orthogonal matrices, we construct $\bm{M} \in \mathbb{R}^{n \times m}$ with entries $M_{ij} \sim \mathcal{N}(0,1)$ i.i.d. and define $\bm{Q}$ from the singular value decomposition of $\bm{G}$, $\bm{Q} = \bm{U}$ where $\bm{G}=\bm{U \Sigma V}^\top$ \citep[Lemma 2.6][]{tulino2004random}.}

We now show that it provides a $1/nm$-approximation guarantee for Problem \eqref{prob:orthogonallyconstrained}. 
\begin{proposition}\label{prop:haar_1}
    Let $\bm{Q} \in \mathbb{R}^{n \times m}$ be distributed uniformly over the set of semi-orthogonal matrices, $\{ \bm{U} \in \mathbb{R}^{n \times m} \: : \: \bm{U}^\top \bm{U} = \bm{I}_m \}$. {Then, for any $\bm{A}\in \mathcal{S}^{nm}_{+}$, we have that:}
        \begin{align*}
        \mathbb{E}[\operatorname{vec}(\bm{Q})^\top \bm{A} \operatorname{vec}(\bm{Q})] 
        \: \leq \: \max_{\bm{U} \in \mathbb{R}^{n \times m}: \bm{U}^\top \bm{U}=\bm{I}_m} \operatorname{vec}(\bm{U})^\top \bm{A} \operatorname{vec}(\bm{U}) 
        \: \leq \: nm \, \mathbb{E}[\operatorname{vec}(\bm{Q})^\top \bm{A} \operatorname{vec}(\bm{Q})].
    \end{align*}
\end{proposition}
\begin{proof}
Because $\bm{Q}$ is feasible for \eqref{prob:orthogonallyconstrained}, we have
\begin{align*}
\operatorname{vec}(\bm{Q})^\top \bm{A} \operatorname{vec}(\bm{Q}) \leq 
\max_{{\bm{U} \in \mathbb{R}^{n \times m} : \bm{U}^\top \bm{U}=\bm{I}_m}} \: \operatorname{vec}(\bm{U})^\top \bm{A} \operatorname{vec}(\bm{U}),
\end{align*}
which leads to the first inequality. 

Furthermore, 
\begin{align*}
\max_{{\bm{U} \in \mathbb{R}^{n \times m} : \bm{U}^\top \bm{U}=\bm{I}_m}} \: \operatorname{vec}(\bm{U})^\top \bm{A} \operatorname{vec}(\bm{U}) \leq \max_{\bm{u} \in \mathbb{R}^{nm} : \| \bm{u} \|{^2} = m} \: \bm{u}^\top \bm{A} \bm{u} = m \lambda_{\max}(\bm{A}) \leq m \operatorname{tr}(\bm{A}).
\end{align*}
To conclude, observe that 
since $\bm{Q}$ is distributed according to the Haar measure, we have $\mathbb{E}[\bm{q}_i \bm{q}_i^\top]=\frac{1}{n}\bm{I}_n$ and $\mathbb{E}[\bm{q}_i \bm{q}_j^\top]=\bm{0}$ for $i \neq j$ \citep[cf.][]{meckes2019random}. Therefore, we have $\mathbb{E}[\operatorname{vec}(\bm{Q})\operatorname{vec}(\bm{Q})^\top ]=\frac{1}{n}\bm{I}_{nm}$ and $\mathbb{E}[\operatorname{vec}(\bm{Q})^\top \bm{A} \operatorname{vec}(\bm{Q})]=\frac{1}{n} \operatorname{tr}(\bm{A})$. 
\end{proof}

\begin{remark} Proposition \ref{prop:haar_1}'s upper bound is tight for uniform sampling. Indeed, if $\bm{A}$ is an identity matrix, then any uniformly sampled $\bm{Q}$ is optimal and the left inequality is tight. Moreover, if $\bm{A}$ is a matrix such that $\bm{A}^{(i,j)}_{i,j}=1$ for $i,j \in [m]$ and $A^{(i,j)}_{\ell_1,\ell_2}=0$ for $\ell_1 \neq i$ or $\ell_2 \neq j$ otherwise, then $\mathrm{tr}(\bm{A})=m$ and an optimal choice of $\bm{U}$ is $\bm{U}_i=\bm{e}_i$, giving both $\max_{\bm{U} \in \mathbb{R}^{n \times m}: \bm{U}^\top \bm{U}=\bm{I}_m} \operatorname{vec}(\bm{U})^\top \bm{A} \operatorname{vec}(\bm{U})=m^2$ and
$nm \, \mathbb{E}[\operatorname{vec}(\bm{Q})^\top \bm{A} \operatorname{vec}(\bm{Q})]=m^2$. 
\end{remark}

\subsection{A Deflation Heuristic} \label{ssec:benchmark.deflation}
We now propose a second baseline inspired by the deflation approach for {Principal Component Analysis \citep[e.g.,][]{mackey2008deflation}}. Namely, we consider the first diagonal block of $\bm{A}$, $\bm{A}^{(1,1)}$, and compute its leading eigenvector. This defines the column $\bm{q}_1$. We then update (or deflate) the other diagonal blocks $\bm{A}^{(j,j)} \leftarrow (\bm{I}_n - \bm{q}_1\bm{q}_1^\top) \bm{A}^{(j,j)} (\bm{I}_n - \bm{q}_1\bm{q}_1^\top)$ and proceed with another block. {We describe the procedure in Algorithm \ref{alg:deflation}, where we assume the diagonal blocks are processed in the natural order, $1,2,\dots,m$.} 
\begin{algorithm}
\caption{A Deflation-Inspired Benchmark 
}
\label{alg:deflation}
\begin{algorithmic}
    \REQUIRE Positive semidefinite matrix $\bm{A} \in \mathcal{S}^{nm}_{+}$
    \STATE Initialize $\bm{Q} = \bm{0}$
    \FOR{$i=1,\dots,m$}
        \STATE Deflate the block $\bm{A}^{(i,i)}$ according to $\bm{q}_1,\dots,\bm{q}_{i-1}$, i.e., compute 
        $$\bm{B} = (\bm{I}_n-\bm{q}_{i-1}\bm{q}_{i-1}^\top) \cdots (\bm{I}_n-\bm{q}_1\bm{q}_1^\top) \bm{A}^{(i,i)} (\bm{I}_n-\bm{q}_1\bm{q}_1^\top) \cdots (\bm{I}_n-\bm{q}_{i-1}\bm{q}_{i-1}^\top)$$
        \STATE Compute $\bm{v}_i \in \arg\max_{\bm{x} \: \| \bm{x} \|_2 = 1} \: \bm{x}^\top \bm{B} \bm{x}$
        \STATE Define $\bm{q}_i = z_i \bm{v}_i$ with $\mathbb{P}(z_i = 1) = 1-\mathbb{P}(z_i = -1) = 1/2$.
    \ENDFOR
    \RETURN Semi-orthogonal matrix $\bm{Q}$ 
\end{algorithmic}
\end{algorithm}
\begin{remark}
    Observe that if $\bm{A}$ is a block diagonal matrix with identical diagonal blocks $\bm{A}^{(i,i)}=\bm{\Sigma}$ and zero off diagonal blocks $\bm{A}^{(i,j)}=\bm{0}$ for $i \neq j$, as in principal component analysis, then the proposed algorithm corresponds to deflation in PCA and is thus exact.
\end{remark}

In practice, we can process the diagonal blocks in any order, with each ordering leading to a different candidate solution. In our implementation, we consider $N$ random permutations of $\{1,\dots,m\}$ and generate $N$ feasible solutions to allow for a fair comparison with our sampling-based approach. 
We show (Proposition \ref{prop:m2eigenapprox}) that it provides a $1/m^2$-factor approximation. 
{We show that Algorithm \ref{alg:deflation}, with a random ordering of the blocks, leads to a $1/m^2$--approximation guarantee.} 
\begin{proposition}\label{prop:m2eigenapprox}
    Let $\bm{Q}$ be generated according to Algorithm \ref{alg:deflation} with a random ordering of the blocks. {Then, for any $\bm{A}\in \mathcal{S}^{nm}_{+}$, we have that:}
    \begin{align*}
        \mathbb{E}\left[\operatorname{vec}(\bm{Q})^\top  \bm{A} \operatorname{vec}(\bm{Q})\right] \leq \max_{\bm{U}\in \mathbb{R}^{n \times m}: \bm{U}^\top \bm{U}=\bm{I}_m}\left[\operatorname{vec}(\bm{U})^\top \bm{A} \operatorname{vec}(\bm{U})\right]\leq m^2 \mathbb{E}\left[\operatorname{vec}(\bm{Q})^\top  \bm{A} \operatorname{vec}(\bm{Q})\right]. 
    \end{align*}
\end{proposition}

\begin{proof}
    First, at iteration $i$, since $\bm{q}_i$ is collinear to the leading eigenvector of the matrix obtained by {projecting} $\bm{A}^{(i,i)}$ onto a space orthogonal to $\bm{q}_1,\dots,\bm{q}_{i-1}$, we have that $\bm{q}_i^\top \bm{q}_j=0$ for each $j < i$ and thus $\bm{Q}$ is semi-orthogonal, so the left inequality holds.

    Second, since $z_i, z_j$ are i.i.d. with mean 0, in expectation, we have that $\mathbb{E}\left[\bm{q}_i^\top \bm{A}^{(i,j)}\bm{q}_j\right]=0$ for  $i \neq j$, and thus the expected objective value attained by $\bm{Q}$ is $$\mathbb{E}[\operatorname{vec}(\bm{Q})^\top \bm{A} \operatorname{vec}(\bm{Q})]=\sum_{i \in [m]}\mathbb{E}[\bm{q}_i^\top \bm{A}^{(i,i)}\bm{q}_i] = \sum_{i \in [m]}\bm{v}_i^\top \bm{A}^{(i,i)}\bm{v}_i. $$ 
    When treating the blocks in the order $1,\dots,m$, we have $\mathbb{E}[\operatorname{vec}(\bm{Q})^\top \bm{A} \operatorname{vec}(\bm{Q})] \geq \bm{v}_1^\top \bm{A}^{(1,1)}\bm{v}_1 = \lambda_{\max}(\bm{A}^{(1,1)})$. By taking the average over random permutations of $\{1,\dots,m\}$, we get $\mathbb{E}[\operatorname{vec}(\bm{Q})^\top \bm{A} \operatorname{vec}(\bm{Q})] \geq \dfrac{1}{m} \sum_{i \in [m]} \lambda_{\max}(\bm{A}^{(i,i)})$.
    
    On the other hand, for any $2 \times 2$ block of $\bm{A}$ {indexed by $(i,j)$ with $i \neq j$},  we have 
    \begin{align*}
    \begin{pmatrix}
        \bm{A}^{(i,i)} & \bm{A}^{(i,j)}\\
        \bm{A}^{(j,i)} & \bm{A}^{(j,j)}
    \end{pmatrix}\succeq \bm{0},
    \end{align*}
    so for any orthogonal matrix $\bm{U}$,
\begin{align*}
\bm{u}_i^\top \bm{A}^{(i,j)}\bm{u}_j + \bm{u}_j^\top \bm{A}^{(j,i)}\bm{u}_i \leq \bm{u}_i^\top \bm{A}^{(i,i)}\bm{u}_i + \bm{u}_j^\top \bm{A}^{(j,j)}\bm{u}_j,
\end{align*}
and
\begin{align*}
\sum_{i,j \in [m]}\bm{u}_i^\top \bm{A}^{(i,j)}\bm{u}_j \leq m \sum_{i \in [m]} \bm{u}_i^\top \bm{A}^{(i,i)}\bm{u}_i \leq m \sum_{i \in [m]} \lambda_{\max}( \bm{A}^{(i,i)}) \leq m^2 \mathbb{E}[\operatorname{vec}(\bm{Q})^\top \bm{A} \operatorname{vec}(\bm{Q})],
\end{align*}
where the last inequality follows from $\mathbb{E}[\operatorname{vec}(\bm{Q})^\top \bm{A} \operatorname{vec}(\bm{Q})] \geq \dfrac{1}{m} \sum_{i \in [m]} \lambda_{\max}(\bm{A}^{(i,i)})$ and concludes the proof.
\end{proof}
\begin{remark}
    Our proof technique shows that  $\mathbb{E}[\operatorname{vec}(\bm{Q})^\top \bm{A} \operatorname{vec}(\bm{Q})] \geq \dfrac{1}{m} \sum_{i \in [m]} \lambda_{\max}(\bm{A}^{(i,i)})$. Thus, Algorithm \ref{alg:deflation} yields a $1/m$-factor approximation when $\bm{A}$ is a block diagonal matrix and a $1/(2m)$-factor approximation when $\bm{A}$ is block diagonally dominant. 
    In general, however, the block diagonal objective bounds the full objective within a factor of $m$, hence the overall $1/m^2$ guarantee. It is also worth noting that the semidefinite inequality
    \begin{align*}
        \begin{pmatrix}
            \bm{A}^{(1,1)} & \bm{A}^{(1,2)} & \ldots \bm{A}^{(1,m)} \\  
            \bm{A}^{(2,1)} & \bm{A}^{(2,2)} & \ldots \bm{A}^{(2,m)} \\  
            \bm{A}^{(3,1)} & \vdots & \ddots \bm{A}^{(3,m)} \\  
            \bm{A}^{(m,1)} & \bm{A}^{(m,2)} & \ldots \bm{A}^{(m,m)} \\  
        \end{pmatrix} \preceq m\begin{pmatrix}
            \bm{A}^{(1,1)} & \bm{0} & \ldots \bm{0} \\  
            \bm{0} & \bm{A}^{(2,2)} & \ldots \bm{0} \\  
            \bm{0} & \vdots & \ddots \bm{0} \\  
            \bm{0} & \bm{0} & \ldots \bm{A}^{(m,m)} \\  
        \end{pmatrix}
    \end{align*}
    which we implicitly prove as part of our approximation guarantee, is actually a special case of the pinching inequality from quantum information theory \citep[see][Lemma II.2]{mosonyi2015quantum}. 
\end{remark}

\section{Numerical Results}\label{sec:numerics}
We numerically evaluate the performance of 
Algorithm~\ref{alg:matrixgwalgorithm2} for semi-orthogonal quadratic optimization problems \eqref{prob:orthogonallyconstrained}. {All experiments are conducted on one Intel(R) Xeon(R) Platinum 8370C 2.80GHz CPU with 128GB of RAM, using the Julia programming language and \verb|Mosek| version 11.0.} For fixed $(n,m)$, we generate a random semidefinite matrix $\bm{A} = \bm{BB}^\top \in \mathcal{S}_+^{nm}$ where the entries of $\bm{B} \in \mathbb{R}^{nm \times 10}$ are {i.i.d. standard normal} random variables. We solve the Shor relaxation \eqref{prob:orth_relax_shor} and sample $N=100$ feasible solutions from Algorithm \ref{alg:matrixgwalgorithm2}. For comparison, we also implement the following benchmarks: 
\begin{itemize}
    \item We sample $N$ solutions uniformly at random (Uniform).
    \item We sample $N$ solutions by applying the earlier deflation heuristic (Deflation).
    \item We follow the heuristic in \citet{burer2023strengthened}, namely, we project the reshaped leading eigenvector of $\bm{W}^\star$ (Burer and Park).
    \item We implement the $1/(2\sqrt{2})$-algorithm described in \citet[][Figure 3]{naor2014efficient}.
\end{itemize}
We consider $n= 100$ and $m \in \{1,2,5,10,15,20,25,30,40, 60, 80, 100\}$. We generate five instances for each $(n,m)$. 

\begin{figure}[h]
    \centering
    \includegraphics[width=\textwidth]{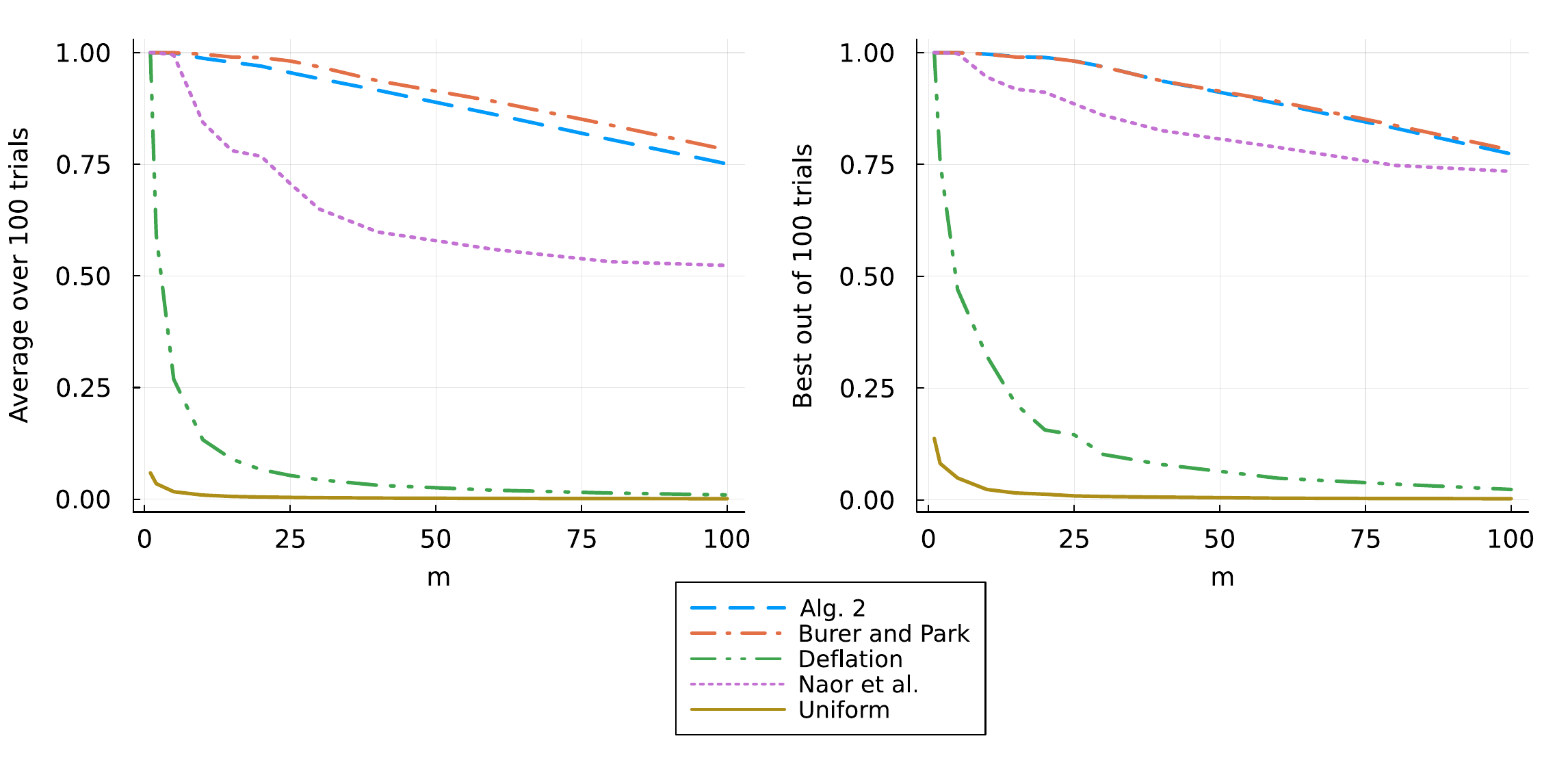}
    \caption{Average approximation ratio $\operatorname{vec}(\bm{Q})^\top  \bm{A} \operatorname{vec}(\bm{Q}) / \langle \bm{A}, \bm{W}^\star \rangle$ (left panel) and approximation ratio $\operatorname{vec}(\bm{Q})^\top  \bm{A} \operatorname{vec}(\bm{Q}) / \langle \bm{A}, \bm{W}^\star \rangle$ of the best solution (right panel) over $N=100$ samples, for different feasibility heuristics. Note that the method of \citet{burer2023strengthened} only returns one solution. For each value of $m$, results are averaged over 5 instances.}
    \label{fig:performance_ratio.panel}
\end{figure}

The left panel of Figure \ref{fig:performance_ratio.panel} compares the average approximation ratio for these {five} algorithms. Confirming our theoretical analysis, we observe that our Algorithm \ref{alg:matrixgwalgorithm2} strongly outperforms Deflation and Uniform---theoretically, Uniform and Deflation achieve a $1/nm$- and $1/m^2$-approximation ratio respectively (Proposition \ref{prop:haar_1} and \ref{prop:m2eigenapprox}). {Despite its mildly stronger worst-case theoretical approximation ratio ($1/(2\sqrt{2})$ vs $1/3$), the algorithm of \citet{naor2014efficient} achieves a worst average empirical performance.}
We also find that the heuristic in \citet{burer2023strengthened} {achieves strong performance, comparable to that of Algorithm \ref{alg:matrixgwalgorithm2}}.

A more relevant performance metric in practice is the performance of the best solution found (out of $N$), rather than the average performance. Regarding the best solution found, the right panel of Figure \ref{fig:performance_ratio.panel} shows that the relative ordering of the methods remains unchanged, although the gap between methods shrinks.   

{To support the above discussion, Table} \ref{tab:orthquad.relaxtime} reports the time required to solve our semidefinite relaxation \eqref{prob:orth_relax_shor} for different values of $m$, using Mosek as the semidefinite optimization solver. Table \ref{tab:orthoquad.sampling} reports the time required by each feasibility heuristic (excluding time to solve the relaxation when needed). {We observe that using the \verb|Mosek| solver, our relaxation requires $1$--$2$ orders of magnitude more computational time to solve than to run benchmark heuristics like deflation or uniform rounding, but facilitates the selection of significantly higher quality feasible solutions.}
\begin{table}[h]
    \centering
    \caption{Computational time (average and standard deviation) for solving \eqref{prob:orth_relax_shor} for $n=100$ and various values of $m$. Results are aggregated over 5 instances.}
    \label{tab:orthquad.relaxtime}
\begin{tabular}{rrr}
$m$ & Average Time (s) & Std Dev (s) \\
\midrule
1 & 2.01 & 0.15 \\
2 & 3.71 & 0.34 \\
5 & 7.35 & 0.63 \\
10 & 19.06 & 1.69 \\
15 & 38.07 & 2.94 \\
20 & 80.33 & 7.34 \\
25 & 140.67 & 21.99 \\
30 & 244.0 & 37.19 \\
40 & 510.24 & 82.15 \\
60 & 1789.24 & 321.35 \\
80 & 3240.3 & 552.06 \\
100 & 7233.63 & 69.91 \\
\end{tabular}
\end{table}

\begin{table}[h]
\centering
\caption{Computational time (average and standard deviation) for different feasibility heuristics for $n=100$ and various values of $m$. For methods that require solving the relaxation \eqref{prob:orth_relax_shor} (Alg. 2, \cite{burer2023strengthened}, \citet{naor2014efficient}), time for solving the relaxation is not included but reported in Table \ref{tab:orthquad.relaxtime}. Results are aggregated over 5 instances.}
\label{tab:orthoquad.sampling}

\begin{tabular}{cccccc}
$m$ & Alg. 2 & Burer and Park & Deflation & Naor et al. & Uniform \\
\hline
2 & 0.01 (0.0) & 0.01 (0.0) & 0.23 (0.09) & 49.25 (0.77) & 0.01 (0.0) \\
5 & 0.07 (0.01) & 0.06 (0.01) & 0.8 (0.45) & 48.27 (0.76) & 0.07 (0.01) \\
10 & 0.34 (0.03) & 0.22 (0.02) & 2.7 (1.76) & 47.43 (0.35) & 0.33 (0.03) \\
15 & 0.6 (0.35) & 0.5 (0.02) & 5.81 (3.98) & 48.13 (0.79) & 0.53 (0.19) \\
20 & 0.99 (0.31) & 0.92 (0.01) & 9.17 (6.44) & 47.34 (0.77) & 0.92 (0.15) \\
25 & 1.44 (0.24) & 1.49 (0.02) & 9.32 (0.86) & 47.91 (0.54) & 1.44 (0.24) \\
30 & 2.59 (0.35) & 2.36 (0.03) & 12.1 (0.53) & 48.31 (1.43) & 2.59 (0.36) \\
40 & 7.0 (0.71) & 5.51 (0.1) & 20.37 (1.03) & 50.92 (1.04) & 7.0 (0.71) \\
60 & 26.98 (1.93) & 17.01 (0.24) & 43.14 (2.38) & 73.5 (3.46) & 26.97 (1.93) \\
80 & 79.83 (23.52) & 37.32 (0.3) & 77.44 (6.37) & 111.68 (3.59) & 79.82 (23.52) \\
100 & 124.65 (3.67) & 66.53 (0.38) & 116.52 (9.99) & 172.1 (7.01) & 124.34 (3.54) \\
\end{tabular}
\end{table}

\section{Conclusion} 
This paper proposes a new technique for relaxing and rounding quadratic optimization problems over semi-orthogonal matrices, which {generalizes} 
the blueprint of the Goemans--Williamson algorithm for BQO. 
{Our algorithm provides a constant-factor approximation} for the orthogonally constrained quadratic optimization problem \eqref{prob:orthogonallyconstrained}, which subsumes the heterogeneous PCA problem \eqref{prob:orthogonallyconstrained_blockcase} among others. Future work could investigate the theoretical analysis of other sampling schemes or extend the approach (namely, solving a Shor relaxation followed by a sampling algorithm) to a broader class of problems, such as low-rank optimization problems.

\putbib[thebib]
\end{bibunit}

\newpage 
\FloatBarrier
\begin{appendices}
\numberwithin{theorem}{section} 
\numberwithin{equation}{section} 
\numberwithin{table}{section} 
\numberwithin{remark}{section} 
\numberwithin{algorithm}{section} 


\setcounter{algorithm}{0} 

\begin{center}
    \bf \Large Improved Approximation Algorithms for Orthogonally Constrained Problems Using Semidefinite Optimization (Appendix)
\end{center}

\begin{bibunit}[plainnat]
  \renewcommand*{\bibnumfmt}[1]{[R#1]}
  \renewcommand*{\citenumfont}[1]{R#1}

\section{Connection with Binary Quadratic Optimization and the Original Goemans--Williamson Algorithm} \label{sec:a.bqo.equiv}
We connect our quadratic semi-orthogonal optimization problem \eqref{prob:orthogonallyconstrained}, its semidefinite relaxation \eqref{prob:orth_relax_shor}, and our rounding algorithm (Algorithm \ref{alg:matrixgwalgorithm2}) to the canonical binary quadratic optimization problem. 

{First, we show the following reduction result between Problems \eqref{prob:maxcut2} and \eqref{prob:orthogonallyconstrained}: 
\begin{proposition} \label{prop:red.preserving}
    Consider an instance of the binary quadratic optimization problem \eqref{prob:maxcut2} with $\bm{Q} \succeq \bm{0}$. We can construct an optimization problem of the form \eqref{prob:orthogonallyconstrained} with $n \geq m$ and such that: 
    \begin{itemize}
        \item any feasible solution to  \eqref{prob:maxcut2} can be converted (in polynomial time) into a feasible solution to \eqref{prob:orthogonallyconstrained} with the same objective value; 
        \item any feasible solution to \eqref{prob:orthogonallyconstrained} can be converted (in polynomial time) into a feasible solution to \eqref{prob:maxcut2} with equal or higher objective value.
    \end{itemize}
\end{proposition}
}
\begin{proof} Consider a binary quadratic optimization problem \eqref{prob:maxcut2}: {
\begin{align*}
    \max_{\bm{z} \in \{-1,1\}^m} \: \bm{z}^\top \bm{Q}\bm{z}.
\end{align*}
Fix an integer $n \geq m$ and denote $\{ \bm{e}_i\}_{i=1,\dots,n}$ the canonical basis of $\mathbb{R}^n$. Define the $(i,j)$ block of the matrix $\bm{A}$ as} $\bm{A}^{(i,j)} := Q_{i,j} \bm{e}_i \bm{e}_j^\top$. In particular, we have that $\bm{Q} \succeq \bm{0} \iff \bm{A} \succeq \bm{0}$. {Consider the corresponding instance of Problem \eqref{prob:orthogonallyconstrained}. 

For any feasible solution to \eqref{prob:maxcut2}, we can construct a feasible solution to \eqref{prob:orthogonallyconstrained} with equal cost. Indeed,} for each $i  \in [m]$, {we can} define $\bm{u}_i := z_i \bm{e}_i$. By construction, we have $\bm{u}_i^\top \bm{u}_j = 0$ if $i \neq j$ and $\bm{u}_i^\top \bm{u}_i = z_i^2 = 1$. With this {construction}, 
\begin{align*}
    \bm{z}^\top \bm{Q}\bm{z} = \sum_{i,j \in [m]} Q_{i,j} z_i z_j = \sum_{i,j \in [m]}  (\bm{u}_i^\top \bm{e}_i) Q_{i,j} (\bm{e}_j^\top \bm{u}_j) = \sum_{i,j \in [m]} \bm{u}_i^\top \bm{A}^{(i,j)} \bm{u}_j.
\end{align*}

{
Alternatively, consider a feasible solution to this instance of \eqref{prob:orthogonallyconstrained}, with objective value 
\begin{align*}
\sum_{i,j \in [m]} \bm{u}_i^\top \bm{A}^{(i,j)} \bm{u}_j = \sum_{i,j \in [m]}  (\bm{u}_i^\top \bm{e}_i) Q_{i,j} (\bm{e}_j^\top \bm{u}_j) = \sum_{i,j \in [m]}  u_{i,i} Q_{i,j}  u_{j,j}.
\end{align*}
We show how to construct a feasible solution $\bm{z} \in \{-1,1\}^m$ to \eqref{prob:maxcut2} with objective greater than or equal to $\sum_{i,j}  u_{i,i} Q_{i,j}  u_{j,j}$. Start from $z_i := u_{i,i} \in [-1,1]$. Fix $z_i$ for all $i>1$. The function $z_1 \in [-1,1] \mapsto \sum_{i,j}  Q_{i,j} z_i z_j$ is convex because $\bm{Q} \succeq \bm{0}$. Thus we can shift $z_1$ to one of the endpoints, $-1$ or $1$, without decreasing the value. Proceeding in this way with the other coordinates $z_2,\dots,z_m$, we construct a solution $\bm{z} \in \{-1,1\}^m$ such that $\bm{z}^\top \bm{Q} \bm{z} \geq \sum_{i,j}  u_{i,i} Q_{i,j}  u_{j,j}$.
}
\end{proof}

{
First of all, Proposition \ref{prop:red.preserving} reduces any BQO problem with a positive semidefinite objective matrix $\bm{Q} \succeq \bm{0}$ to a problem of the form \eqref{prob:orthogonallyconstrained} with equal objective value. We note that our procedure to construct a feasible solution to BQO from a solution to \eqref{prob:orthogonallyconstrained} is analogous to that of \citet{naor2014efficient} and leverages the convexity of the objective. In particular, the standard (nonnegative weight) Max-Cut can be formulated as a BQO problem of this class, with $\bm{Q}$ being the (weighted) Laplacian of the graph. So Proposition \ref{prop:red.preserving} shows that Problem \eqref{prob:orthogonallyconstrained} is NP-hard. \citet[][Theorem 3.1]{lai2025stiefel} provide an alternative proof of this result. They also use a reduction from Max-Cut. However, their BQO formulation of Max-Cut is different, and the corresponding matrix $\bm{Q}$ is not PSD so they need to use a different decoding scheme. 

Second, the reduction in Proposition \ref{prop:red.preserving} is  approximation-preserving: Any polynomial-time 
$\alpha$-approximation (in objective value) for Problem \eqref{prob:orthogonallyconstrained} would immediately yield an 
$\alpha$-approximation for Max-Cut.  Therefore, 
the inapproximability threshold of ${2/\pi}+\varepsilon$ for BQO \citep[][Theorem I.3]{briet2015tight}
transfers to Problem~\eqref{prob:orthogonallyconstrained}.}

{In the rest of this section, we show how our semidefinite relaxation and our {relax-and-project} algorithm would work in the special case where Problem~\eqref{prob:orthogonallyconstrained} encodes an instance of \eqref{prob:maxcut2}.

First, in this case, let us observe that optimal solutions to \eqref{prob:orthogonallyconstrained} can be found among solutions of the form $\bm{u}_i = z_i \bm{e}_i$ with $z_i \in \{-1,1\}$. In other words, although we do not impose the constraint that each vector $\bm{u}_i$  should be colinear to $\bm{e}_i$, optimality naturally enforces this constraint. 
\begin{proof}
For any feasible solution to \eqref{prob:orthogonallyconstrained}, its objective value is $\sum_{i,j} \bm{u}_i^\top \bm{A}^{(i,j)} \bm{u}_j = \sum_{i,j}  u_{i,i} Q_{i,j}  u_{j,j}$. Denoting $z_i := u_{i,i} \in [-1,1]$, we have that 
\begin{align*}
    \eqref{prob:orthogonallyconstrained} \leq \max_{\bm{z} \in [-1,1]^m} \: \bm{z}^\top \bm{Q}\bm{z} = \max_{\bm{z} \in \{-1,1\}^m} \: \bm{z}^\top \bm{Q}\bm{z},
\end{align*}
where the last equality follows from the fact that $\bm{Q} \succeq \bm{0}$. 
Conversely, for each $\bm{z} \in \{-1,1\}^m$ the matrix defined as $\bm{u}_i = z_i \bm{e}_i$ is feasible for \eqref{prob:orthogonallyconstrained} and achieves an objective value of $\bm{z}^\top \bm{Q}\bm{z}$.
\end{proof}
}

{We now generalize this observation to the semidefinite relaxation of \eqref{prob:orthogonallyconstrained} and show that, without loss of optimality, the semidefinite variable $\bm{W}$ is of the form $\bm{W}^{(i,j)} = Z_{i,j} \bm{e}_i \bm{e}_j^\top$ for some matrix $\bm{Z} \in \mathcal{S}^m_+$.
\begin{proof}
For any feasible solution to the semidefinite relaxation \eqref{prob:orth_relax_shor}, its objective value is 
$$\langle \bm{W}, \bm{A} \rangle = \sum_{i,j \in [m]} \langle \bm{W}^{(i,j)} , \bm{A}^{(i,j)}  \rangle = \sum_{i,j \in [m]} Q_{i,j} \bm{e}_i^\top \bm{W}^{(i,j)} \bm{e}_j, = \langle \bm{Z}, \bm{Q} \rangle,$$
with $Z_{i,j} := \bm{e}_i^\top \bm{W}^{(i,j)} \bm{e}_j$. The constraint $\bm{W} \succeq \bm{0}$ implies $\bm{Z} \succeq \bm{0}$ and the constraints $\operatorname{tr}(\bm{W}^{(j,j)}) = 1,\ j\in [m]$ imply $Z_{j,j} \leq 1, j\in [m]$. For any $j \in [m]$, the function $Z_{j,j} \in [0,1] \mapsto \langle \bm{Z}, \bm{Q} \rangle$ is linear with slope $Q_{j,j} \geq 0$ so, for any feasible $\bm{Z}$, setting $Z_{j,j}$ to 1 cannot decrease the objective value (and does not break feasibility). Consequently, at optimality, we can assume that $Z_{j,j} =W^{(j,j)}_{j,j} = 1$. However, $\operatorname{tr}(\bm{W}^{(j,j)}) = 1$. Thus, all other diagonal coefficients of the block $\bm{W}^{(j,j)}$ are equal to 0. Recall that, for a positive semidefinite matrix, if a diagonal coefficient is equal to 0, then its entire column/row has to be equal to 0. Because $\bm{W}^{(j,j)} \succeq \bm{0}$, it means that 
$\bm{W}^{(j,j)} = Z_{j,j} \bm{e}_j \bm{e}_j^\top$. Regarding the off-diagonal blocks of $\bm{W}$, for $(i,j)$ with $i \neq j$, because $\bm{W} \succeq \bm{0}$, the columns $j'$ of $\bm{W}^{(i,j)}$ with $j' \neq j$ are equal to 0, i.e., only the $j$th column of $\bm{W}^{(i,j)}$ can be nonzero. Similarly, all the rows $i'$ with $i' \neq i$ of $\bm{W}^{(i,j)}$ are equal to 0. All in all, we have that the blocks of $\bm{W}$ are of the form $\bm{W}^{(i,j)} = Z_{i,j} \bm{e}_i\bm{e}_j^\top$.
\end{proof}}
{Considering a solution to the semidefinite relaxation} of the form $\bm{W}^{(i,j)} = Z_{i,j} \bm{e}_i \bm{e}_j^\top$. The objective of \eqref{prob:orth_relax_shor} can thus be written as 
\begin{align*}
    \langle \bm{A}, \bm{W} \rangle = \sum_{i,j} \langle \bm{A}^{(i,j)}, \bm{W}^{(i,j)} \rangle = \sum_{i,j} Q_{i,j} Z_{i,j},
\end{align*}
and the constraints on the matrix $\bm{W}$ are equivalent to: 
\begin{align*}
    & \bm{W} \succeq \bm{0}: & \bm{Z} \succeq \bm{0}, \\
    &\operatorname{tr}(\bm{W}^{(j,j')}) = \delta_{j,j'}: & Z_{j,j} = 1,\\
    &\sum_{i \in [m]} \bm{W}^{(i,i)} \preceq \bm{I}_n: & Z_{i,i} \leq 1, \ \forall i \in [m].
\end{align*}
So, we recover the semidefinite relaxation of BQO, \eqref{prob:maxcut3}, exactly. 

Consider a solution to the semidefinite relaxation \eqref{prob:orth_relax_shor}, $\bm{W}^\star$, and $\bm{Z}^\star$ such that $\bm{W}^{\star(i,j)} = Z^\star_{i,j} \bm{e}_i \bm{e}_j^\top$. By sampling $\operatorname{vec}(\bm{G}) \sim \mathcal{N}(\bm{0}, \bm{W}^\star)$, the sparsity pattern of $\bm{W}^\star$ implies that each column of $\bm{G}$, $\bm{g}_i$, is of the form $\bm{g}_i = y_i \bm{e}_i$, with $\bm{y} \sim \mathcal{N}(\bm{0},\bm{Z}^\star)$. In this case, the matrix $\bm{G}$ is diagonal and its SVD can be written 
\begin{align*}
    \bm{G} = \bm{U \Sigma V}^\top := {\begin{pmatrix} \bm{I}_m \\ \bm{0}_{(n-m) \times m} \end{pmatrix}} \: \begin{pmatrix} |y_1| & & \\ & \ddots & & \\ & & |y_m| \end{pmatrix} \begin{pmatrix} \operatorname{sign}(y_1) & & \\ & \ddots & \\ & & \operatorname{sign}(y_m) \end{pmatrix}.
\end{align*}
We then generate { $\bm{Q} = \bm{U V}^\top$ so} 
each column of $\bm{Q}$, $\bm{q}_i$, can be expressed as $\bm{q}_i = \operatorname{sign}(y_i) \bm{e}_i$, i.e., $\hat{z}_i = \operatorname{sign}(y_i)$, which is precisely the original Goemans--Williamson algorithm (Algorithm \ref{alg:gwmethod}). 

{\blue
\section{Proof of Proposition \ref{prop:ce.13}} \label{ssec:a.proof.ce13}
The construction of our instance relies on a group of orthogonal matrices known as the real Clifford group \citep[see, e.g.,][]{hashagen2018real}. We first show an intermediate result.

\begin{lemma}\label{lemma:r.clifford} Consider an integer $m = 2^q$ for some $q \geq 1$. Denoting $n := 1 + 2^{q^2 + q+ 2}(2^q -1) \prod_{j\in [q-1]}(4^j - 1)$, there exist $n-1$ matrices $\bm{R}_i \in \mathcal{S}^m$ of the form $\bm{R}_i = \bm{U}_i \bm{D} \bm{U}_i^\top$ with $\bm{U}_i \in \mathbb{R}^{m \times m}$ orthogonal and $\bm{D} = \operatorname{Diag}(\bm{I}_{m/2}, -\bm{I}_{m/2})$ such that 
\begin{align} \label{eqn:clifford.lambda}
   \lambda_{\max}\left([\tfrac{1}{m} \operatorname{tr}(\bm{R}_i^\top \bm{R}_j]_{i,j}\right) < \dfrac{n-1}{m-1},
\end{align}
and, for any unit-norm vector $\bm{u} \in \mathbb{R}^m$,
\begin{align} \label{eqn:clifford.avg}
    \dfrac{1}{n-1} \sum_{i\in [n-1]} \bm{R}_i \bm{uu}^\top \bm{R}_i^\top = \dfrac{m}{(m-1)(m+2)} \bm{I}_m + \dfrac{m-2}{(m-1)(m+2)} \bm{uu}^\top.
\end{align}
\end{lemma}
\begin{proof} The proof relies on the real Clifford group, $\mathcal{H}$, a finite subset of the set of $m \times m$ orthogonal matrices \citep[see, e.g.,][]{hashagen2018real}. Denote $n-1 := | \mathcal{H} | = 2^{q^2 + q+ 2}(2^q -1) \prod_{j\in [q-1]}(4^j - 1)$ \citep[see][Equation 6.3.4]{nebe2006self} and enumerate $\mathcal{H} = \{ \bm{H}_1,\dots, \bm{H}_{n-1}\}$. The real Clifford group is a finite set of orthogonal matrices satisfying (among others) a property referred to as the exact twirling identity \citep[see, e.g., Equation (50) in][]{hashagen2018real}: For any matrix $\bm{X} \in \mathbb{R}^{m^2 \times m^2}$, 
\begin{align*}
    \dfrac{1}{|\mathcal{H}|} \sum_{\bm{H} \in \mathcal{H}} (\bm{H} \otimes \bm{H}) \bm{X} (\bm{H} \otimes \bm{H})^\top = \mathbb{E}_{\bm{U}}[(\bm{U} \otimes \bm{U}) \bm{X} (\bm{U} \otimes \bm{U})^\top],
\end{align*}
where $\bm{U}$ is distributed according to the uniform (Haar) measure over the set of orthogonal matrices. In other words, the average of $(\bm{H} \otimes \bm{H}) \bm{X} (\bm{H} \otimes \bm{H})^\top$ over the finite set of matrices $\mathcal{H}$ equals the average over $\bm{H}$ distributed according to the Haar measure. \citet[][Section IV.B]{audenaert2002asymptotic} give an explicit expression for the Haar-average, namely
\begin{align*} 
    \bm{T}(\bm{X}) &:= \mathbb{E}_{\bm{U}}[(\bm{U} \otimes \bm{U}) \bm{X} (\bm{U} \otimes \bm{U})^\top] \\ &= \operatorname{tr}(\bm{P}_\Omega \bm{X}) \bm{P}_\Omega + \dfrac{2}{m(m-1)} \operatorname{tr}(\bm{P}_A \bm{X}) \bm{P}_A + \dfrac{2}{(m+2)(m-1)} [\operatorname{tr}((\bm{P}_S - \bm{P}_\Omega) \bm{X})] (\bm{P}_S - \bm{P}_\Omega)
\end{align*}
where 
\begin{itemize}
    \item $\bm{P}_\Omega$ is the orthogonal projection onto the span of $\operatorname{vec}(\bm{I}_m)$, i.e., $\bm{P}_\Omega = \operatorname{vec}(\bm{I}_m) \operatorname{vec}(\bm{I}_m)^\top / m$ 
    \item $\bm{P}_A$ is the orthogonal projection onto the span of $\{ \operatorname{vec}(\bm{A}) \: : \: \bm{A}^\top = -\bm{A} \}$,
    \item $\bm{P}_S$ is the orthogonal projection onto the span of $\{ \operatorname{vec}(\bm{S}) \: : \: \bm{S}^\top = \bm{S} \}$.
\end{itemize}
In particular, if $\bm{K} \in \mathbb{R}^{m^2 \times m^2}$ denotes the commutation matrix \citep[see, e.g.,][]{magnus1979commutation}, i.e., the matrix $\bm{K} $ such that, for any $\bm{X} \in \mathbb{R}^{m \times m}$, $\operatorname{vec}(\bm{X}^\top) = \bm{K} \operatorname{vec}(\bm{X}) $, then $\bm{P}_A = (\bm{I}_{m^2} - \bm{K})/2$ and $\bm{P}_S = (\bm{I}_{m^2} + \bm{K})/2$.

For $i=1,\dots,n-1$, we define $\bm{R}_i = \bm{H}_i \bm{D} \bm{H}_i^\top \in \mathcal{S}^m$. Let us show that the matrices $\bm{R}_i$ constructed in this manner satisfy the desired properties.

To prove the first property, let us denote $\bm{v}_i = \operatorname{vec}(\bm{R}_i) = (\bm{H}_i \otimes \bm{H}_i) \operatorname{vec}(\bm{D})$ and $C_{i,j} := \bm{v}_i^\top \bm{v}_j /m$. With this notation, 
\begin{align*}
    \lambda_{\max}(\bm{C}) &= \dfrac{1}{m} \lambda_{\max}(\bm{V}^\top \bm{V}) = \dfrac{1}{m} \lambda_{\max}(\bm{V}\bm{V}^\top ) \\ &= \dfrac{n-1}{m} \lambda_{\max}\left(\dfrac{1}{|\mathcal{H}|} \sum_{\bm{H} \in \mathcal{H}} (\bm{H} \otimes \bm{H})\operatorname{vec}(\bm{D})\operatorname{vec}(\bm{D})^\top (\bm{H} \otimes \bm{H}) \right).
\end{align*}
The matrix on the right-hand side is exactly an average over the real Clifford group, hence it is equal to $\bm{T}(\operatorname{vec}(\bm{D})\operatorname{vec}(\bm{D})^\top)$. 
We compute 
\begin{align*}
\operatorname{tr}(\bm{P}_\Omega\operatorname{vec}(\bm{D})\operatorname{vec}(\bm{D})^\top) &= \dfrac{1}{m} (\operatorname{vec}(\bm{D})^\top \operatorname{vec}(\bm{I}_m))^2 = \dfrac{1}{m} (\operatorname{tr}(\bm{D}))^2 = 0, \\
\operatorname{tr}(\bm{P}_A \operatorname{vec}(\bm{D})\operatorname{vec}(\bm{D})^\top) &= \operatorname{tr}(\bm{0}) = 0, \\
\operatorname{tr}(\bm{P}_S \operatorname{vec}(\bm{D})\operatorname{vec}(\bm{D})^\top) &= \operatorname{vec}(\bm{D})^\top\operatorname{vec}(\bm{D}) = m,
\end{align*}
so 
\begin{align*} 
    \bm{T}(\operatorname{vec}(\bm{D})\operatorname{vec}(\bm{D})^\top) &= \dfrac{2m}{(m+2)(m-1)}  (\bm{P}_S - \bm{P}_\Omega).
\end{align*}
Consequently, 
\begin{align*}
    \lambda_{\max}(\bm{C}) \leq\dfrac{n-1}{m} \dfrac{2m}{(m+2)(m-1)} = \dfrac{n-1}{m-1} \dfrac{2}{m+2} < \dfrac{n-1}{m-1},
\end{align*}
for $m \geq 2$ (or, equivalently, $q \geq 1$). 

For the second property, let us consider the vectorized version of the matrix identity
\begin{align*}
     \dfrac{1}{n-1} \sum_{i\in [n-1]} \operatorname{vec}(\bm{R}_i \bm{uu}^\top \bm{R}_i^\top) &= \dfrac{1}{n-1} \sum_{i\in [n-1]} (\bm{R}_i \otimes \bm{R}_i) \operatorname{vec}(\bm{uu}^\top) 
     \\ &= \left(
     \dfrac{1}{|\mathcal{H}|} \sum_{\bm{H} \in \mathcal{H}}(\bm{H} \otimes \bm{H}) (\bm{D} \otimes \bm{D}) (\bm{H} \otimes \bm{H})\right) \operatorname{vec}(\bm{uu}^\top) 
     \\ &= \bm{T}(\bm{D} \otimes \bm{D}) \operatorname{vec}(\bm{uu}^\top). 
\end{align*}
We have
\begin{align*}
\operatorname{tr}(\bm{P}_\Omega(\bm{D} \otimes \bm{D})) &= \dfrac{1}{m} \operatorname{vec}(\bm{D}\bm{I}_m\bm{D})^\top \operatorname{vec}(\bm{I}_m) = 1.
\end{align*}
For $\bm{P}_A(\bm{D} \otimes \bm{D})$ and $\bm{P}_S(\bm{D} \otimes \bm{D})$, we use the fact that the commutation matrix satisfies $\operatorname{tr}(\bm{K}(\bm{A} \otimes \bm{B})) = \operatorname{tr}(\bm{A} \bm{B})$ for any matrices $\bm{A},\bm{B} \in \mathbb{R}^{m\times m}$ \citep[Theorem 3.1(xiii) in][]{magnus1979commutation} to write
\begin{align*}
\operatorname{tr}(\bm{P}_A(\bm{D} \otimes \bm{D})) &= \dfrac{1}{2} \left( \operatorname{tr}(\bm{D} \otimes \bm{D}) - \operatorname{tr}(\bm{K}(\bm{D} \otimes \bm{D})) \right) = - \dfrac{1}{2} m,\\
\operatorname{tr}(\bm{P}_S(\bm{D} \otimes \bm{D})) &= \dfrac{1}{2} \left( \operatorname{tr}(\bm{D} \otimes \bm{D}) + \operatorname{tr}(\bm{K}(\bm{D} \otimes \bm{D})) \right) = \dfrac{1}{2} m.
\end{align*}
As a result, 
\begin{align*} 
    \bm{T}(\bm{D} \otimes \bm{D}) &= \bm{P}_\Omega - \dfrac{1}{m-1} \bm{P}_A + \dfrac{m-2}{(m+2)(m-1)} (\bm{P}_S - \bm{P}_\Omega).
\end{align*}

The matrix $\bm{uu}^\top$ is symmetric, so $\bm{P}_A  \operatorname{vec}(\bm{uu}^\top) = \bm{0}$ and $\bm{P}_S  \operatorname{vec}(\bm{uu}^\top) = \operatorname{vec}(\bm{uu}^\top)$. Furthermore, 
\begin{align*}
    \bm{P}_{\Omega} \operatorname{vec}(\bm{uu}^\top) = \dfrac{1}{m} \operatorname{vec}(\bm{I}_m)\operatorname{vec}(\bm{I}_m)^\top  \operatorname{vec}(\bm{uu}^\top) = \dfrac{1}{m} \operatorname{tr}(\bm{uu}^\top) \operatorname{vec}(\bm{I}_m) = \dfrac{1}{m} \operatorname{vec}(\bm{I}_m).
\end{align*}
Therefore, we proved that 
\begin{align*}
    \operatorname{vec}\left( \dfrac{1}{n-1} \sum_{i\in [n-1]} \bm{R}_i \bm{uu}^\top \bm{R}_i^\top \right) &= \bm{T}(\bm{D} \otimes \bm{D}) \operatorname{vec}(\bm{uu}^\top) \\
    &= \dfrac{1}{m} \left( 1 - \dfrac{m-2}{(m+2)(m-1)}\right) \operatorname{vec}(\bm{I}_m) +  \dfrac{m-2}{(m+2)(m-1)} \operatorname{vec}(\bm{uu}^\top) \\
        &=  \dfrac{m}{(m+2)(m-1)} \operatorname{vec}(\bm{I}_m) +  \dfrac{m-2}{(m+2)(m-1)} \operatorname{vec}(\bm{uu}^\top)
\end{align*}
which concludes the proof.
\end{proof}

With this intermediate result, we can turn to the proof of Proposition~\ref{prop:ce.13}.
\begin{proof}[Proof of Proposition~\ref{prop:ce.13}] Fix $m = 2^q$ with $q \geq 1$. By Lemma \ref{lemma:r.clifford}, there exists $n \geq m$ and symmetric matrices $\bm{R}_2,\dots,\bm{R}_n$ of the form $\bm{R}_i = \bm{U}_i \bm{D} \bm{U}_i^\top$ with $\bm{U}_i \in \mathbb{R}^{m \times m}$ orthogonal and $\bm{D} = \operatorname{Diag}(\bm{I}_{m/2}, -\bm{I}_{m/2})$ satisfying \eqref{eqn:clifford.lambda} and \eqref{eqn:clifford.avg}.

For this $n$, let $\bm{e}_1\in\mathbb{R}^n$ denote the first standard basis vector in $\mathbb{R}^n$ and define $\bm{A} := \bm{I}_m\otimes (\bm{e}_1 \bm{e}_1^\top) = \operatorname{Diag}(\bm{e}_1 \bm{e}_1^\top, \dots, \bm{e}_1 \bm{e}_1^\top) \in\ \mathcal{S}^{nm}_+$.

Given that $\bm{A}$ is block-diagonal with identical diagonal blocks, Problem \eqref{prob:orthogonallyconstrained} is equivalent to finding the top-$m$ eigenvectors of $\bm{e}_1 \bm{e}_1^\top$. The optimal objective value is 1 (e.g., achieved by the semi-orthogonal matrix $\bm{U} = [\bm{e}_1,\dots,\bm{e}_m]$). Similarly, the objective value of the semidefinite relaxation \eqref{prob:orth_relax_shor} is 1: for any feasible $\bm{W} \in \mathcal{S}_+^{nm}$, we have 
$$ \langle \bm{A}, \bm{W} \rangle = \sum_{j=§}^m \langle \bm{e}_1 \bm{e}_1^\top, \bm{W^{(j,j)}} \rangle =  \langle \bm{e}_1 \bm{e}_1^\top, \sum_{j=§}^m \bm{W^{(j,j)}} \rangle \leq \langle \bm{e}_1 \bm{e}_1^\top, \bm{I}_n \rangle = 1,$$
and the reverse inequality holds because it is a relaxation.
For any matrix $\bm{Q} \in \mathbb{R}^{n \times m}$, its objective value is
$\operatorname{vec}(\bm{Q})^\top \bm{A} \operatorname{vec}(\bm{Q}) = \sum_{j\in [m]} (\bm{q}_j^\top \bm{e}_1)^2 = \| \bm{e}_1^\top \bm{Q}\|_2^2$, i.e., it is equal to the squared norm of the first row of $\bm{Q}$. 

We construct the matrix $\bm{G}$ for Algorithm \ref{alg:matrixgwalgorithm2} row by row: we sample $\bm{x} \sim \mathcal{N}(\bm{0},\bm{I}_m)$ and define 
\begin{align*}
\bm{G} = \dfrac{1}{\sqrt{m}} \begin{pmatrix} \bm{x}^\top \\ \rho \, \bm{x}^\top \bm{R}_2 \\ \vdots \\ 
\rho \, \bm{x}^\top \bm{R}_n
\end{pmatrix}
\quad \mbox{ or equivalently } \quad 
\bm{G}^\top  = \dfrac{1}{\sqrt{m}} \begin{bmatrix} \bm{x}  \ | \ \rho \bm{R}_2^\top \bm{x} \ | \ \cdots \ | \ \rho \bm{R}_n^\top\bm{x}
\end{bmatrix},
\end{align*}
with $\rho = \sqrt{(m-1)/(n-1)} \leq 1$.
The vector $\operatorname{vec}(\bm{G})$ is normally distributed, $\mathbb{E}[\operatorname{vec}(\bm{G})] = \bm{0}$, so we have implicitly constructed $\bm{W} := \mathbb{E}[\operatorname{vec}(\bm{G})\operatorname{vec}(\bm{G})^\top] \succeq \bm{0}$. In the remainder of this proof, we show that this matrix $\bm{W}$ is feasible and optimal for the semidefinite relaxation, and compute $\operatorname{vec}(\bm{Q})^\top \bm{A} \operatorname{vec}(\bm{Q})$ analytically. 

First, let us show that $\bm{W}$ is feasible for the SDP relaxation . The matrix $\bm{W}$ is feasible if and only if the corresponding random matrix $\bm{G}$ satisfies $\mathbb{E}[\bm{G}^\top \bm{G} ] = \bm{I}_m$ and $\mathbb{E}[\bm{G} \bm{G}^\top ] \preceq \bm{I}_n$. In our case, we have 
\begin{align*}
\bm{G}^\top \bm{G} = \sum_{i\in [n]} \bm{G}^\top \bm{e}_i \bm{e}_i^\top \bm{G} = \dfrac{1}{m} \bm{xx}^\top + \dfrac{\rho^2}{m}  \sum_{i=2}^n \bm{R}_i^\top \bm{xx}^\top \bm{R}_i.
\end{align*}
Taking expectations and using the facts that $\mathbb{E}[\bm{xx}^\top] = \bm{I}_m$ and $\bm{R}_i^\top \bm{R}_i = \bm{U}_i^\top \bm{D}^2 \bm{U}_i = \bm{I}_m$, we have
$$\mathbb{E}[\bm{G}^\top \bm{G}] = \dfrac{1}{m} \left(1 + (n-1) \rho^2\right) \bm{I}_m = \bm{I}_m.$$
In addition, 
\begin{align*}
\mathbb{E}\left[ \bm{G} \bm{G}^\top \right]= \begin{pmatrix}
    1 & \bm{0} \\ 
    \bm{0} & \rho^2 \bm{C}
\end{pmatrix} \mbox{ with } C_{i,j} := \dfrac{1}{m} \operatorname{tr}(\bm{R}_i^\top \bm{R}_{j} ).
\end{align*}
By construction (Equation \eqref{eqn:clifford.lambda} in Lemma \ref{lemma:r.clifford}), $\lambda_{\max}(\bm{C}) < \dfrac{n-1}{m-1} = \rho^{-2}$ so $\mathbb{E}\left[ \bm{G} \bm{G}^\top \right] \preceq \bm{I}_n$. Hence, $\bm{W}$ is a feasible solution to the semidefinite relaxation \eqref{prob:orth_relax_shor}.

Second, we have that 
$\langle \bm{A}, \bm{W} \rangle = \langle \bm{e}_1 \bm{e}_1^\top, \sum_{j=§}^m \bm{W^{(j,j)}} \rangle  = \langle \bm{e}_1 \bm{e}_1^\top, \mathbb{E}[\bm{G}\bm{G}^\top ] \rangle = 1$. So $\bm{W}$ is optimal.

Finally, let us evaluate the performance of $\bm{Q}$. By definition, $\bm{Q} = \bm{G (\bm{G}^\top \bm{G})^{-1/2} }$ so 
\begin{align*}
   \| \bm{e}_1^\top \bm{Q}\|_2^2 =  \bm{e}_1^\top \bm{Q}\bm{Q}^\top \bm{e}_1 = \bm{e}_1^\top  \bm{G} (\bm{G}^\top \bm{G})^{-1} \bm{G}^\top \bm{e}_1 = \dfrac{1}{m} \bm{x}^\top (\bm{G}^\top \bm{G})^{-1}  \bm{x}.
\end{align*}
Furthermore, 
\begin{align*}
\bm{G}^\top \bm{G} = \dfrac{1}{m} \bm{xx}^\top + \dfrac{\rho^2}{m}  \sum_{i=2}^n \bm{R}_i^\top \bm{xx}^\top \bm{R}_i.
\end{align*}
Denoting $\bm{u} = \bm{x} / \| \bm{x} \|$, we have 
\begin{align*}
\bm{G}^\top \bm{G} &= \| \bm{x} \|^2 \left(\dfrac{1}{m} \bm{uu}^\top + \dfrac{\rho^2}{m}  \sum_{i=2}^n \bm{R}_i^\top \bm{uu}^\top \bm{R}_i \right), \\
  \| \bm{e}_1^\top \bm{Q}\|_2^2 &= \bm{u}^\top\left( \bm{uu}^\top + {\rho^2} \sum_{i=2}^n \bm{R}_i^\top \bm{uu}^\top \bm{R}_i \right)^{-1}  \bm{u}.
\end{align*}
By construction (Equation \eqref{eqn:clifford.avg} in Lemma \ref{lemma:r.clifford}), we have, for any unit-norm vector $\bm{u}$,
\begin{align*}
    \dfrac{1}{n-1} \sum_{i=2}^n \bm{R}_i^\top \bm{uu}^\top \bm{R}_i =  \dfrac{m}{(m-1)(m+2)} \bm{I}_m +  \dfrac{m-2}{(m-1)(m+2)} \bm{uu}^\top.
\end{align*}
So, 
\begin{align*}
\bm{uu}^\top + {\rho^2} \sum_{i=2}^n \bm{R}_i^\top \bm{uu}^\top \bm{R}_i
\quad =  \dfrac{m}{m+2} \bm{I}_m +  \dfrac{2m}{m+2} \bm{uu}^\top.
\end{align*}
In particular, we get that $\bm{u}$ is an eigenvector of the matrix above with eigenvalue $3m/(m+2)$ so we have
\begin{align*}
\bm{u}^\top \left( \bm{uu}^\top + {\rho^2} \sum_{i=2}^n \bm{R}_i^\top \bm{uu}^\top \bm{R}_i \right)^{-1} \bm{u}
\quad =  \dfrac{m+2}{3m}.
\end{align*}
In other words, $\operatorname{vec}(\bm{Q})^\top \bm{A}\operatorname{vec}(\bm{Q}) = \dfrac{m+2}{3m}$ almost surely. 
\end{proof}
}

\section{Technical Appendix to Section \ref{sec:proof}} \label{sec:a.proof}

\subsection{Proof of Lemma \ref{lemma:moments}} \label{ssec:a.moments}
\begin{proof} As described in Section \ref{ssec:gw_matrix.sampling}, in our implementation of Algorithm \ref{alg:matrixgwalgorithm2}, we sample $\operatorname{vec}(\bm{G}) \sim \mathcal{N}(\bm{0}_{nm}, \bm{W}^\star)$ as $\operatorname{vec}(\bm{G}) = \sum_{k \in [r]} \operatorname{vec}(\bm{B}_k) z_k$ 
with $\bm{z} \sim \mathcal{N}(\bm{0}_r, \bm{I}_r)$ and $\bm{W}^\star = \sum_{k \in [r]} \operatorname{vec}(\bm{B}_k)\operatorname{vec}(\bm{B}_k)^\top$ a Cholesky decomposition of $\bm{W}^\star$. This construction interprets $\bm{G}$ as a matrix series, $\bm{G} = \sum_{k \in [r]} \bm{B}_k z_k$, as studied in the statistics literature \citep[see, e.g.,][]{tropp2015introduction}.

{Let us denote 
\begin{align*}
    \bm{V}_1 := \quad \sum_{k \in [r]} \bm{B}_k\bm{B}_k^\top \mbox{ and } \quad 
    \bm{V}_2 := \quad \sum_{k \in [r]} \bm{B}_k^\top \bm{B}_k ,
\end{align*}
so $\mathbb{E}[\bm{GG}^\top ] = \bm{V}_1$ and  $\mathbb{E}[\bm{G}^\top \bm{G} ] = \bm{V}_2$.}
Noting that $\bm{W}^{(i,j)} = \sum_{k \in [r]} \bm{B}_k \bm{e}_i \bm{e}_j^\top \bm{B}_k^\top$, we have
\begin{align*}
    & \mathbb{E}\left[\bm{G}^\top \bm{G} \right]_{i,j} = \left( \sum_{k \in [r]} \bm{B}_k^\top \bm{B}_k \right)_{i,j} \quad = \sum_{k \in [r]} \bm{e}_i^\top \bm{B}_k^\top \bm{B}_k \bm{e}_j \quad = \operatorname{tr}(\bm{W}^{(i,j)}), \\
    \mbox{and } & \mathbb{E}\left[\bm{G} \bm{G}^\top \right] = \sum_{k \in [r]} \bm{B}_k \bm{B}_k^\top \quad = \sum_{k \in [r]} \sum_{i \in [m]} \bm{B}_k \bm{e}_i \bm{e}_i^\top \bm{B}_k^\top \quad = \sum_{i \in [m]} \bm{W}^{(i,i)}.
\end{align*}
The fact that $\bm{W}$ satisfies the constraints in \eqref{prob:orth_relax_shor} yields $\bm{V}_1 \preceq \bm{I}_n$ and $\bm{V}_2 = \bm{I}_m$, {as claimed}. 

For the fourth moments, we define the following \emph{Hermitian} Gaussian series
\begin{align*}
    \bm{Y} := \sum_{k \in [r]} z_k \bm{A}_k \ \mbox{ with } \ \bm{A}_k := \begin{pmatrix}
        \bm{0} & \bm{B}_k \\ \bm{B}_k^\top & \bm{0}
    \end{pmatrix}.
\end{align*}
By construction,
\begin{align*}
    \sum_{k \in [r]} \bm{A}_k^2 
    = \begin{pmatrix}
        \sum_k \bm{B}_k\bm{B}_k^\top & \bm{0} \\ \bm{0} &  \sum_k \bm{B}_k^\top\bm{B}_k \end{pmatrix} 
    = \begin{pmatrix} \bm{V}_1 & \bm{0} \\ \bm{0} &  \bm{V}_2 \end{pmatrix} 
\quad \mbox{ and } \quad 
\bm{Y}^4 =
\begin{pmatrix}
(\bm{GG}^\top)^2 & 0\\
0 & (\bm{G}^\top \bm{G})^2
\end{pmatrix}.
\end{align*}
Hence it suffices to upper bound $\mathbb E[\bm{Y}^4]$ and then read
the diagonal blocks.

Since $\bm{z}\sim\mathcal N(\bm{0},\bm{I}_r)$, Isserlis' theorem gives
$\mathbb E[z_a z_b z_c z_d] =\delta_{ab}\delta_{cd}+\delta_{ac}\delta_{bd}+\delta_{ad}\delta_{bc}$ \citep{isserlis1918formula} .
Expanding $\bm{Y}^4$and taking expectation yields
\begin{align*}
\mathbb E[\bm{Y}^4]
&=\sum_{a,b,c,d \in [r]} \mathbb E[z_a z_b z_c z_d] \bm{A}_a \bm{A}_b \bm{A}_c \bm{A}_d\\
&=\sum_{a,c \in [r]} \bm{A}_a^2 \bm{A}_c^2
+\sum_{a,b \in [r]} \bm{A}_a \bm{A}_b \bm{A}_a \bm{A}_b
+\sum_{a,b \in [r]} \bm{A}_a \bm{A}_b^2 \bm{A}_a\\
&=\Big(\sum_{k \in [r]} \bm{A}_k^2\Big)^2
+\sum_{a,b \in [r]} \bm{A}_a \bm{A}_b \bm{A}_a \bm{A}_b
+\sum_{a,b \in [r]} \bm{A}_a \bm{A}_b^2 \bm{A}_a.
\end{align*}

Fix two indices $a,b$. By expanding the inequality $(\bm{A}_a \bm{A}_b - \bm{A}_b \bm{A}_a)^\top(\bm{A}_a \bm{A}_b - \bm{A}_b \bm{A}_a)\succeq 0$ and using the fact that $\bm{A}_a,\bm{A}_b$ are symmetric, we have
\begin{align*}
    \bm{A}_a \bm{A}_b \bm{A}_a \bm{A}_b + \bm{A}_b \bm{A}_a \bm{A}_b \bm{A}_a
\preceq
\bm{A}_a \bm{A}_b^2 \bm{A}_a + \bm{A}_b \bm{A}_a^2 \bm{A}_b.
\end{align*}
Summing over all $(a,b)$ 
yields
\begin{align*}
\sum_{a,b \in [r]} \bm{A}_a \bm{A}_b \bm{A}_a \bm{A}_b \preceq \sum_{a,b \in [r]} \bm{A}_a \bm{A}_b^2 \bm{A}_a. 
\end{align*}
Therefore,
\[
\mathbb E[\bm{Y}^4]
\preceq
\Big(\sum_{k \in [r]} \bm{A}_k^2\Big)^2
+2\sum_{a,b \in [r]} \bm{A}_a \bm{A}_b^2 \bm{A}_a
=
\Big(\sum_{a \in [r]} \bm{A}_a^2\Big)^2
+2\sum_{a \in [r]} \bm{A}_a\Big(\sum_{b=1}^r \bm{A}_b^2\Big)\bm{A}_a.
\]
Reading the first diagonal block, we have 
\[
\mathbb{E}[(\bm{GG}^\top)^2] \preceq \bm{V}_1^2 + 2 \sum_{k \in [r]} \bm{B}_k \bm{V}_2\bm{B}_k^\top \preceq \bm{V}_1^2 + 2 \lambda_{\max}(\bm{V}_2) \bm{V}_1 \preceq 3 \bm{I}_n.
\]
Similarly for the second diagonal block, 
\[
\mathbb{E}[(\bm{G}^\top\bm{G})^2] \preceq \bm{V}_2^2 + 2 \sum_{k \in [r]} \bm{B}_k^\top  \bm{V}_1 \bm{B}_k \preceq \bm{V}_2^2 + 2 \lambda_{\max}(\bm{V}_1) \bm{V}_2 \preceq 3 \bm{I}_m,
\]
which concludes the proof.
\end{proof}

{\blue
\subsection{Proof of Lemma \ref{lemma:matrixgaussian.12}}\label{ssec:a.12}
\begin{proof} Denote $\bm{N}$ a square root of $\bm{M}$. By Cauchy-Schwarz,
\begin{align*}
   \operatorname{tr}(\bm{M}) = \operatorname{tr}(\bm{M}\mathbb{E}[\bm{S}]) = \mathbb{E}[\operatorname{tr}(\bm{N} \bm{S}^{1/4} \bm{S}^{3/4} \bm{N})] \leq \sqrt{\mathbb{E}[\operatorname{tr}(\bm{N} \bm{S}^{1/2} \bm{N})] } \sqrt{\mathbb{E}[\operatorname{tr}(\bm{N} \bm{S}^{3/2} \bm{N})] }.
\end{align*}
Applying Cauchy-Schwarz again, we have
\begin{align*}
   \mathbb{E}[\operatorname{tr}(\bm{N} \bm{S}^{3/2} \bm{N})] = \mathbb{E}[\operatorname{tr}(\bm{N} \bm{S} \bm{S}^{1/2} \bm{N})] \leq  \sqrt{\mathbb{E}[\operatorname{tr}(\bm{N} \bm{S}^2 \bm{N})]} \sqrt{\mathbb{E}[\operatorname{tr}(\bm{N} \bm{S} \bm{N})]} \leq \sqrt{c} \operatorname{tr}(\bm{M}).
\end{align*}
Combining these inequalities together, we have
\begin{align*}
   \operatorname{tr}(\bm{M})^2 \leq \mathbb{E}[\operatorname{tr}(\bm{N} \bm{S}^{1/2} \bm{N})] \, \sqrt{c} \operatorname{tr}(\bm{M}).
\end{align*}
Rearranging the terms concludes the proof.
\end{proof}
}

{
\section{A tighter dimension-dependent approximation ratio}}
\begin{theorem}\label{thm:chisquared.light}
{Assume that $\bm{A} \succeq \bm{0}$. Let $\bm{W}^\star$ denote an optimal solution to the semidefinite relaxation \eqref{prob:orth_relax_shor}. The random matrix $\bm{Q} \in \mathbb{R}^{n \times m}$ generated by Algorithm~\ref{alg:matrixgwalgorithm2} satisfies the inequality
\begin{align*}
\mathbb{E}\left[\operatorname{vec}(\bm{Q})^\top \bm{A} \operatorname{vec}(\bm{Q}) \right] & \geq  \beta_{n,m} \langle \bm{A}, \bm{W}^\star \rangle.
\end{align*}
}
with 
\begin{align*}
\beta_{n,m} := \min_{ \lambda \in [0,1]} \: \int_{0}^\infty \left(1 + 2 t \, m\dfrac{1-\lambda}{nm-1} \right)^{-(nm-1)/2} \, (1 + 2 t m \lambda)^{-3/2}dt.
\end{align*}

In particular, the constant $\beta_{n,m}$ satisfies the following properties:
\begin{itemize}
    \item[(a)] For any integer $m$, $\beta_{n,m}$ is non-increasing in $n$. For any integer $n$, $\beta_{n,m}$ is non-increasing in $m$. 
    \item[(b)] For any integer $m$, we have $\beta_{n,m} \rightarrow \beta_{\infty,m}$ as $n \rightarrow \infty$ with 
\begin{align*}
\beta_{\infty,m} &:= \min_{ \lambda \in [0,1]} \: \int_{0}^\infty e^{- t \, m (1-\lambda)} \, (1 + 2 t m \lambda)^{-3/2} dt 
\end{align*}
\item[(c)] For $m=1$, $\beta_{n,1}$ is optimal, i.e., there exists a covariance matrix $\bm{W}^\star$ satisfying $\mathbb{E}\left[\operatorname{vec}(\bm{Q})^\top \bm{A} \operatorname{vec}(\bm{Q}) \right] =\beta_{n,1} \langle \bm{A}, \bm{W}^\star \rangle$.
\end{itemize}
\end{theorem}
\begin{proof}
{From the proof of Proposition \ref{prop:2pim}, we have
\begin{align*}
    \mathbb{E}[\operatorname{vec}(\bm{Q})^\top \bm{A} \operatorname{vec}(\bm{Q})] \geq \langle \bm{A}, \mathbb{E}\left[ \dfrac{\operatorname{vec}(\bm{G}) \operatorname{vec}(\bm{G})^\top}{\|\bm{G}\|_F^2} \right] \rangle.
\end{align*}
}

We use the {integral representation $1/x^2 = \int_{t \geq 0} e^{-t x^2} dt$} to obtain 
\begin{align*}
    \mathbb{E}\left[\frac{\operatorname{vec}(\bm{G})\operatorname{vec}(\bm{G})^\top}{\|\bm{G}\|_F^2}\right]
    &=\int_{0}^\infty \mathbb{E}[\operatorname{vec}(\bm{G})\operatorname{vec}(\bm{G})^\top e^{-t \| \bm{G}\|_F^2}]dt.
\end{align*}
Denote $r = \operatorname{rank}(\bm{W}^\star) \leq nm$ and consider an eigenvalue decomposition of $\bm{W}^\star$, $\bm{W}^\star = \bm{H \Lambda H}^\top$. We have $\operatorname{vec}(\bm{G}) = \bm{H \Lambda}^{1/2}\bm{z}$, with $\bm{z} \sim \mathcal{N}(\bm{0}_r,\bm{I}_r)$ and thus
\begin{align*}
\mathbb{E}[ \operatorname{vec}(\bm{G})\operatorname{vec}(\bm{G})^\top\exp(-t \Vert \bm{G}\Vert_F^2)] 
&= \bm{H\Lambda}^{1/2} \mathbb{E}[ \bm{zz}^\top\exp(-t \bm{z}^\top \bm{\Lambda} \bm{z})]  \bm{\Lambda}^{1/2} \bm{H}^\top.
\end{align*}
Furthermore, {for $\bm{z} \sim \mathcal{N}(\bm{0}_r,\bm{I}_r)$, we can compute the expectation analytically:} 
\begin{align*}
\mathbb{E}[ \bm{zz}^\top\exp(-t \bm{z}^\top \bm{\Lambda} \bm{z})]  
= \dfrac{1}{\sqrt{\operatorname{det}(\bm{I}_{r} + 2t \bm{\Lambda} )}} \, ( \bm{I}_r + 2t \bm{\Lambda})^{-1}.
\end{align*}
So, the integral is lower bounded by
\begin{align*}
     \bm{H\Lambda}^{1/2} \bm{B}  \bm{\Lambda}^{1/2} \bm{H}^\top \mbox{ with } \bm{B} := \int_{0}^\infty \dfrac{1}{\sqrt{\operatorname{det}(\bm{I}_{r} + 2t \bm{\Lambda} )}} \, ( \bm{I}_r + 2t \bm{\Lambda})^{-1} dt.
\end{align*}
To conclude, we show that there exists a scalar $\beta > 0$ such that $\bm{B} \succeq \beta \bm{I}_r$. 

To find such $\beta$, observe that $\bm{B}$ is a diagonal matrix.
Hence, it is sufficient to find a lower bound on its diagonal entries. 
Given the constraints on $\bm{W}^\star$, the eigenvalues $\bm{\Lambda}$ must satisfy: $\Lambda_i \geq 0$ (from $\bm{W}^\star \succeq 0$), $\sum_{i=1}^r \Lambda_i = m$ (from $\operatorname{tr}(\bm{W}^\star) = \sum_{i\in [m]} \operatorname{tr}(\bm{W}^{\star (i,i)}) =  m$). 
Hence, we can take
\begin{align*}
\beta = \min_{\bm{\Lambda} \in [0,m]^r \: : \: \sum_{i=1}^r \Lambda_i = m} \: \int_{0}^\infty \prod_{i=1}^r \left(1 + 2 t \Lambda_{i'} \right)^{-1/2} \, (1 + 2 t \Lambda_{1})^{-1}dt.
\end{align*}
The function $\bm{\Lambda} \mapsto \int_{0}^\infty \prod_{i=1}^r \left(1 + 2 t \Lambda_{i'} \right)^{-1/2} \, (1 + 2 t \Lambda_{1})^{-1}dt$ is convex, and 
invariant under any permutation of the $\Lambda_i$, $i>1$. 
So, by Jensen's inequality, we can restrict our attention to minimizers of the form $\Lambda_{1}=\lambda$, $\Lambda_{i} = \dfrac{m-\lambda}{r-1}, i>1$: 
\begin{align*}
\beta = \min_{ \lambda \in [0,m]} \: \int_{0}^\infty \left(1 + 2 t \, \dfrac{m-\lambda}{r-1} \right)^{-(r-1)/2} \, (1 + 2 t \lambda)^{-3/2}dt.
\end{align*}

For a fixed value of $(t,\lambda)$, the integrand is decreasing in $r \leq nm$, so 
\begin{align*}
\beta \geq \beta_{n,m} & :=  \min_{ \lambda \in [0,m]} \: \int_{0}^\infty \left(1 + 2 t \, \dfrac{m-\lambda}{nm-1} \right)^{-(nm-1)/2} \, (1 + 2 t \lambda)^{-3/2}dt. 
\end{align*}
The change of variable $\lambda \gets \lambda/m$ {leads to the desired expression for $\beta_{n,m}$}. 

Claim (a) follows from the fact that for any $\lambda \in [0,1]$, the integrand is decreasing in $n$ for $m$ fixed, and decreasing in $m$ for $n$ fixed. 

Claim (b) follows from the fact that the sequence $(1+x/k)^{-k}$ converges monotonically to $e^{-x}$ combined with the dominated convergence theorem. 

Claim (c) follows from the fact that, for $m = 1$, $\bm{Q} = \bm{G} / \| \bm{G} \|_F$ so our analysis is exact.
\end{proof}
\begin{remark}\label{rk:tighter.beta} We observe that the bound $\beta_{n,m}$ is obtained by looking at the worst-case instance over all covariance matrices $\bm{W}^\star$. In particular, we can obtain tighter values of $\beta_{n,m}$ by allowing for dependency on the rank of $\bm{W}^\star$, $r$, instead of the ambient dimension $n$. For instance, if $r=1$, we can get 
\begin{align*} 
\beta = \min_{ \lambda \in [0,m]} \: \int_{0}^\infty \, (1 + 2 t \lambda)^{-3/2}dt = \min_{ \lambda \in [0,m]} \: \dfrac{1}{\lambda} = \dfrac{1}{m}.
\end{align*}
Alternatively, by the Barvinok-Pataki bound, we know there exists some optimal solution $\bm{W}^\star$ with rank at most $n+m$ and we could use this bound to refine our constant. 
Furthermore, if $\bm{W}^\star$ has additional structure (e.g., $\bm{W}^\star$ is block diagonal), we can derive additional constraints on the eigenvalues $\bm{\Lambda}$, hence tighter constants $\beta_{n,m}$.
\end{remark}
Compared with Proposition~\ref{prop:2pim}, the value of 
Theorem \ref{thm:chisquared.light} is primarily computational. By solving numerically the one-dimensional minimization problem in $\lambda$, it provides tighter estimates of the performance of our algorithm, especially for small values of $m$, as reported in Table \ref{tab:guarantee.value}. While the guarantee from Proposition~\ref{prop:2pim} is independent of $n$, the constant $\beta_{n,m}$ in Theorem~\ref{thm:chisquared.light} is monotonically decreasing with $n$, obtaining stronger approximation ratios for finite values of $n$. 
However, we should acknowledge that $\beta_{n,m}$ is worse than $1/3$ for $m > 2$ (actually, Remark \ref{rk:tighter.beta} identifies a class of matrices $\bm{W}^\star$ for which $\beta_{n,m} \leq 1/m$).
\begin{table}
    \centering
    \caption{Values of the approximation factor from Proposition~\ref{prop:2pim} and Theorem~\ref{thm:chisquared.light} for some values of $n$ and $m$, with $m \leq n$.} \label{tab:guarantee.value}
\begin{tabular}{r r *{4}{r}}
\toprule
    & \multirow{2}{*}{Proposition \ref{prop:2pim}} & \multicolumn{4}{c}{$\beta_{n,m}$ (Theorem \ref{thm:chisquared.light})} \\
\cmidrule(l){3-6}
    $m$ &  & $n=5$ & $n=10$ & $n=15$ & $n=\infty$ \\
\midrule
1  & 0.636620 & 0.735264 & 0.706972 & 0.697920 & 0.680415 \\
2  & 0.318310 & 0.353486 & 0.346734 & 0.344533 & 0.340208 \\
3  & 0.212207 & 0.232640 & 0.229689 & 0.228720 & 0.226805 \\
4  & 0.159155 & 0.173367 & 0.171721 & 0.171179 & 0.170104 \\
5  & 0.127324 & 0.138164 & 0.137116 & 0.136770 & 0.136083 \\
10 & 0.063662 & \multicolumn{1}{c}{—} & 0.068299 & 0.068213 & 0.068042 \\
15 & 0.042441 & \multicolumn{1}{c}{—} & \multicolumn{1}{c}{—} & 0.045437 & 0.045361 \\
\bottomrule
\end{tabular}
\end{table}


\FloatBarrier

\putbib[thebib]
\end{bibunit}
\end{appendices}

\end{document}